\documentclass[leqno]{amsart}
\usepackage{amsmath,amsthm,amssymb,mathrsfs}
\usepackage[all,cmtip]{xy}

\numberwithin{equation}{subsection}
\newcommand{\setno}{
  \setcounter{equation}{\value{subsubsection}}
  \stepcounter{subsubsection}
}

\newcommand{\N}{\mathbb{N}}
\newcommand{\Q}{\mathbb{Q}}
\newcommand{\R}{\mathbb{R}}
\newcommand{\Z}{\mathbb{Z}}
\newcommand{\BF}{\mathbb{F}}

\newcommand{\BP}{\mathbb{P}}

\newcommand{\D}{\mathcal{D}}
\newcommand{\E}{\mathcal{E}}
\newcommand{\F}{\mathcal{F}}
\newcommand{\K}{\mathcal{K}}
\newcommand{\M}{\mathcal{M}}
\newcommand{\CB}{\mathcal{B}}
\newcommand{\CC}{\mathcal{C}}
\newcommand{\CL}{\mathcal{L}}
\newcommand{\CN}{\mathcal{N}}

\newcommand{\CQ}{\mathcal{Q}}

\newcommand{\FA}{\mathfrak{A}}
\newcommand{\FD}{\mathfrak{D}}
\newcommand{\FE}{\mathfrak{E}}
\newcommand{\FF}{\mathfrak{F}}
\newcommand{\FG}{\mathfrak{G}}
\newcommand{\FK}{\mathfrak{K}}

\newcommand{\FM}{\mathfrak{M}}
\newcommand{\FQ}{\mathfrak{Q}}
\newcommand{\FZ}{\mathfrak{Z}}

\newcommand{\RHom}{\mathcal{RH}om}
\newcommand{\pH}{{}^pH}
\newcommand{\ptau}{{}^p\tau}

\newcommand{\mof}{\text{-mof}}
\newcommand{\pmof}{\text{-pmof}}
\newcommand{\pt}{\text{pt}}

\newcommand{\hI}{\hat{I}}
\newcommand{\hH}{\hat{H}}
\newcommand{\hi}{\hat\imath}
\newcommand{\hj}{\hat\jmath}
\newcommand{\tE}{\tilde{E}}
\newcommand{\tZ}{\tilde{Z}}
\newcommand{\tp}{\tilde{p}}
\newcommand{\TE}{\tilde\E}
\newcommand{\TF}{\tilde\F}
\newcommand{\TN}{\tilde\CN}
\newcommand{\Gn}{G_\nu}
\newcommand{\Gnp}{G_{\nu'}}
\newcommand{\Gnn}{G_{\nu\nu'}}
\newcommand{\nO}{_{\nu,\Omega}}
\newcommand{\npO}{_{\nu',\Omega}}
\newcommand{\Oi}{_{\Omega,i}}
\newcommand{\nOi}{_{\nu,\Omega,i}}
\newcommand{\npOi}{_{\nu',\Omega,i}}
\newcommand{\Ok}{_{\Omega,k}}
\newcommand{\nOk}{_{\nu,\Omega,k}}
\newcommand{\npOk}{_{\nu',\Omega,k}}
\newcommand{\dx}{\dot{x}}
\newcommand{\dy}{\dot{y}}
\newcommand{\FRes}{\mathfrak{Res}}
\newcommand{\vm}{{\vec\mu}}
\newcommand{\vn}{{\vec\nu}}
\newcommand{\vo}{{\vec\omega}}
\newcommand{\bG}{\mathbf{G}}
\newcommand{\tG}{\mathbf{\tilde G}}
\newcommand{\von}{_{\vo,\nu}}
\newcommand{\vonp}{_{\vo,\nu'}}
\newcommand{\Ao}{\FA^\bullet_\vo}

 \DeclareMathOperator{\Ext}{Ext}
 \DeclareMathOperator{\Hom}{Hom}
 \DeclareMathOperator{\Ker}{Ker}
 \DeclareMathOperator{\Coker}{Coker}
 \DeclareMathOperator{\Img}{Img}
 \DeclareMathOperator{\Id}{Id}

 \DeclareMathOperator{\Tr}{Tr}

\newtheorem{thm}[subsubsection]{Theorem}
\newtheorem{prop}[subsubsection]{Proposition}
\newtheorem{lem}[subsubsection]{Lemma}
\newtheorem{cor}[subsubsection]{Corollary}

\theoremstyle{remark}
\newtheorem{rem}[subsubsection]{Remark}
\theoremstyle{definition}

\newcommand{\thmref}[1]{Theorem \ref{#1}}
\newcommand{\propref}[1]{Proposition \ref{#1}}
\newcommand{\lemref}[1]{Lemma \ref{#1}}

\newcommand{\remref}[1]{Remark \ref{#1}}
\newcommand{\secref}[1]{Section \ref{#1}}

\begin{document}

\title[categorification of integrable representations]
{Categorification of integrable representations of quantum groups}
\author[H. Zheng]{Hao Zheng}
\address{Department of Mathematics, Zhongshan University, Guangzhou, Guangdong 510275, China}
\email{zhenghao@mail.sysu.edu.cn}
\thanks{The author is supported in part by an NSFC grant}
\maketitle

\begin{abstract}
We categorify the highest weight integrable representations and
their tensor products of a symmetric quantum Kac-Moody algebra. As
byproducts, we obtain a geometric realization of Lusztig's canonical
bases of these representations as well as a new positivity result.
The main ingredient in the underlying geometric construction is a
class of micro-local perverse sheaves on quiver varieties.
\end{abstract}

\addtocontents{toc}{\protect\setcounter{tocdepth}{1}}

\section*{Introduction}

\subsection{}\label{sec:intro:result}

The purpose of this paper is to categorify integrable
representations of quantum groups. More precisely, given a symmetric
quantum Kac-Moody algebra $U$ and a tensor product $\Lambda$ of
highest weight integrable $U$-modules, we construct the followings.
(See \ref{sec:rep:abel}.)

\begin{enumerate}
\item
An abelian category $\CC$ whose Grothendieck group $\bG(\CC)$ has a
natural structure of free $\Z[q,q^{-1}]$-module.
\item
A set of exact endofunctors of $\CC$ each for a generator of $U$.
\item
A set of functor isomorphisms each for a defining relation of $U$,
so that $\bG(\CC)\otimes\Q(q)$ is endowed with a structure of
$U$-module.
\item
An isomorphism of $U$-modules $\bG(\CC)\otimes\Q(q) \to \Lambda$
such that the indecomposable projectives of $\CC$ give rise to the
canonical basis of $\Lambda$.
\end{enumerate}

Among others, we also obtain a positivity result on canonical basis.
See \thmref{thm:rep:basis}(3) and \remref{rem:rep:main}(3).

\subsection{}

For the simplest quantum group $U_q(sl_2)$, categorification of its
representations has been extensively studied in, for instance,
\cite{CK07,FKS06,La08,Zh07} from various points of view. A
categorification of representations of $U_q(sl_n)$ was obtained in
\cite{Su07} (see also \cite{Ch06} for a treatment of level two
representations). But results for other quantum groups are vacant
from the literature (a good reference on this direction is
\cite{KL08}).

The present work follows the lines of \cite{Zh07} and is intimately
related to the works such as \cite{KSa97, Lu90, Lu91, Ma03, Na01,
Sai02, Sav05} on the interplay between quiver varieties and quantum
groups. Some of the relations will be mentioned in the paper.
Especially, many techniques we use here are borrowed from
\cite{Lu91, Zh07}. Knowledge on them could be of much help for
understanding this paper.

We remark that our categorification can be strengthened in a
straightforward manner to yield a set of functor isomorphisms in
bijection to the multiplication rules of a modified quantum
enveloping algebra as was done for the modified enveloping algebra
of $sl_2$ in \cite{BFK99}. It is however a challenging task to work
out a categorification using higher categories as \cite{CR08,La08}.

\subsection{}

The main contribution of this paper may be summarized as the
construction of suitable categories for categorification task. Below
are some underlying ideas hopefully useful for other problems
concerning quantum algebras.

First of all, a common feature arising from the papers \cite{Lu91,
Zh07} is that one has on hand a certain kind of quiver varieties,
each admitting a structure of cotangent bundle, say $T^*M$. Then one
uses categories of sheaves on $M$, rather than on the total
varieties $T^*M$, to fulfil desired goals.

This is essentially a generalization of the fundamental principle
arising at the very beginning of quantum theory in physics: the
classical phase space describing the motion of a particle is the
symplectic space $\R^3\times\R^3$ formed by the position space and
the momentum space, while the quantum phase space is constituted by
functions on a half part of the total classical phase space, either
the position space or the momentum space.

The above observation led us to the attempt of similar treatment for
the quiver varieties (Nakajima's quiver varieties \cite{Na94,Na98})
involved in our problem. Unfortunately, these quiver varieties in
general do not admit structures of cotangent bundles. So, another
idea came into our picture: we find open immersions of these quiver
varieties into cotangent bundles, say $T^*M$, and extract our
desired categories out by localizing categories of sheaves on $M$.

This is analogous to the definition of the ring of functions on a
quasi-affine variety which is a localization of that on the ambient
affine variety.

In fact, such application of localization on categories of sheaves
is not new. It has already been adopted in the construction of
stacks of micro-local perverse sheaves \cite{GMV05,Wa04}. In this
sense, our categorification is based on a class of micro-local
perverse sheaves on quiver varieties.

Due to some technical reasons, for example characteristic variety
has not been defined for an $l$-adic sheaf at present, realization
of the above ideas is not straightforward in this paper. We will
return to this point in \ref{sec:cat_fun:micro}.

\subsection{}

The paper is organized as follows. In \secref{sec:pre} we introduce
various notions and important facts that will be used in this paper,
as well as claim some elementary consequences without proof.

In \secref{sec:cat_fun}, we first construct a triangulated category
$\FD$, which carries all the information of the highest weight
integrable representations and their tensor products of a quantum
group $U$. Then we define a set of endofunctors of $\FD$ and
establish a set of functor isomorphisms categorifying $U$. The main
theorem is stated in \ref{sec:cat_fun:iso}.

In \secref{sec:rep}, we extract subcategories from $\FD$ each for a
tensor product of highest weight integrable $U$-modules. The abelian
categorification claimed in \ref{sec:intro:result} is fulfilled in
the final subsection. It is helpful to have a glace at the main
results stated in \ref{sec:rep:main} before going into specific
details.


\setcounter{tocdepth}{2}

\tableofcontents

\addtocontents{toc}{\protect\setcounter{tocdepth}{2}}

\section{Preliminaries}\label{sec:pre}

\subsection{Quantum groups}\label{sec:pre:u}

\subsubsection{}

Assume a finite graph without circle edges is given. Let $I$ be the
set of vertices and let $H$ be the set of pairs consisting of an
edge and an orientation of it. The pair $(I,H)$ is a {\em quiver} in
the sense that it also defines a directed graph. The elements of $H$
are referred to as {\em arrows}.

We associate to the finite graph a symmetric generalized Cartan
matrix $(a_{ij})_{i,j\in I}$ by setting $a_{ii}=2$ and $-a_{ij}=$
the number of edges joining $i,j$ if $i\ne j$.

\subsubsection{}

The quantum group $U$ associated to the generalized Cartan matrix
$(a_{ij})$ is the $\Q(q)$-algebra defined by the generators $K_i,
K_i^{-1}, E_i, F_i$, $i\in I$ and the relations
\begin{equation*}
\begin{split}
  & K_iK_i^{-1}=K_i^{-1}K_i=1, \quad K_iK_j=K_jK_i, \\
  & K_iE_j=q^{a_{ij}}E_jK_i, \\
  & K_iF_j=q^{-a_{ij}}F_jK_i, \\
  & E_iF_j-F_jE_i=\delta_{ij}\frac{K_i-K_i^{-1}}{q-q^{-1}}, \\
  & \sum_{m=0}^{1-a_{ij}} (-1)^m E_i^{(m)} E_j E_i^{(1-a_{ij}-m)} = 0, \quad i\ne j, \\
  & \sum_{m=0}^{1-a_{ij}} (-1)^m F_i^{(m)} F_j F_i^{(1-a_{ij}-m)} = 0,  \quad i\ne j,
\end{split}
\end{equation*}
where
\begin{equation*}
  E_i^{(n)} := \frac{E_i^n}{[n]_q!}, \quad F_i^{(n)} := \frac{F_i^n}{[n]_q!},
\end{equation*}
$$[n]_q:=\frac{q^n-q^{-n}}{q-q^{-1}}, \quad
  [n]_q!:=[1]_q[2]_q\cdots[n]_q.
$$

It is a Hopf algebra with comultiplication
\begin{equation*}
\begin{split}
  & \Delta K_i = K_i \otimes K_i, \\
  & \Delta E_i = E_i \otimes 1 + K_i \otimes E_i, \\
  & \Delta F_i = F_i \otimes K_i^{-1} + 1 \otimes F_i,
\end{split}
\end{equation*}
counit
\begin{equation*}
  \varepsilon(K_i) = 1, \quad
  \varepsilon(E_i) = \varepsilon(F_i) = 0,
\end{equation*}
and antipode
\begin{equation*}
  S(K_i) = K_i^{-1}, \quad
  S(E_i) = -K_i^{-1}E_i, \quad
  S(F_i) = -F_iK_i.
\end{equation*}

\begin{rem}
The quantum group $U$ is usually defined a little larger when the
generalized Cartan matrix $(a_{ij})$ is not of finite type (cf.
\cite{Lu93}). But this makes no difference for our discussion on
highest weight integrable representations. So we adopt the present
definition.
\end{rem}

\subsubsection{}

Let $\varrho: U\to U^{op}$ be the algebra isomorphism defined on the
generators by
\begin{equation*}
  \varrho(K_i)=K_i, \quad \varrho(E_i)=qK_iF_i, \quad \varrho(F_i)=qK_i^{-1}E_i.
\end{equation*}
By a {\em contravariant form} on a $U$-module $M$ we mean a
symmetric bilinear form
\begin{equation*}
  (\,,\,): M\times M \to \Q(q)
\end{equation*}
satisfying
\begin{equation*}
  (xu,w)=(u,\varrho(x)w) \quad \text{for $x\in U$, $u,w\in M$}.
\end{equation*}

The isomorphism $\varrho$ is compatible with the comultiplication of
$U$:
$$(\varrho\otimes\varrho)\Delta(x)=\Delta\varrho(x) \quad \text{for $x\in U$}.$$
Hence contravariant forms on $U$-modules $M_1,M_2$ automatically
give rise to a contravariant form on the $U$-module $M_1 \otimes
M_2$
\begin{equation*}
  (u_1\otimes u_2,w_1\otimes w_2) := (u_1,w_1)(u_2,w_2) \quad
  \text{for $u_1,w_1\in M_1$, $u_2,w_2\in M_2$}.
\end{equation*}
We always assume tensor product modules are endowed with
contravariant form in this way.

\subsubsection{}

There is a $\Q$-algebra involution $\bar{\;}: U \to U$ determined by
\begin{equation*}
  \bar{q}=q^{-1}, \quad
  \bar{K_i}=K_i^{-1}, \quad
  \bar{E_i}=E_i, \quad
  \bar{F_i}=F_i.
\end{equation*}

\subsubsection{}\label{sec:pre:u:module}

The highest weight integrable representations of $U$ are described
as follows.

Given $\omega=\sum_i{\omega_i}\cdot i\in\N[I]$, we have a {\em Verma
module} $M(\omega)$ which by definition is the quotient of $U$ by
the left ideal generated by $E_i, K_i-q^{\omega_i}$.

Let $\eta_\omega\in M(\omega)$ be represented by the unit of $U$.
There is a unique contravariant form on $M(\omega)$ satisfying
$(\eta_\omega,\eta_\omega)=1$. Its radical $R(\omega)$ is the
largest nontrivial submodule of $M(\omega)$. The quotient
$\Lambda(\omega)=M(\omega)/R(\omega)$ is a highest weight integrable
$U$-module, and all highest weight integrable $U$-modules are
obtained in this way.

Also let $\eta_\omega\in \Lambda(\omega)$ be represented by the unit
of $U$. We will assume $M(\omega)$, $\Lambda(\omega)$ are both
endowed with the unique contravariant form satisfying
$(\eta_\omega,\eta_\omega)=1$.

For a sequence $\vo = (\omega^1,\omega^2,\dots,\omega^t)$ of
elements in $\N[I]$, we set
$$\Lambda(\vo):=\Lambda(\omega^1)\otimes\Lambda(\omega^2)\otimes\dots\otimes\Lambda(\omega^t).$$
In particular, $\Lambda(\emptyset)$ is the trivial $U$-module
$\Q(q)$.

\subsubsection{}\label{sec:pre:u:verma}

The Verma module $M(\omega)$ is universal in the following sense.

Let $M,N$ be $U$-modules and let $\varphi: M\to N$ be a linear
homomorphism such that $E_i\varphi(u)=\varphi(E_iu)$ and
$K_i\varphi(u)=q^{\omega_i}\varphi(K_iu)$ for $i\in I$ and $u\in M$.
Then there exists a unique $U$-module homomorphism $\tilde\varphi:
M\otimes M(\omega) \to N$ such that
$\tilde\varphi(u\otimes\eta_\omega)=\varphi(u)$, $u\in M$.

Moreover, if $\varphi$ preserves some preassigned contravariant
forms on $M,N$, then so does $\tilde\varphi$.

\subsection{Perverse sheaves on algebraic stacks}\label{sec:pre:perverse}

The theory of perverse sheaves on schemes \cite{BBD82} has a natural
version for algebraic stacks. We refer the readers to \cite{Jo93}
for a systematic treatment.

Algebraic stacks arising from this paper are actually very
restrictive: each is either (i) a (quasi-projective) algebraic
variety or (ii) a quotient stack $[X/G]$ where $X$ is an algebraic
variety and $G$ is a product of general linear groups acting on it.
In the latter case, the derived category of sheaves on $[X/G]$ is
nothing new but the $G$-equivariant derived category described in
\cite{BL94}.

Adopting the language of algebraic stack allows us to treat usual
and equivariant derived categories in a uniform way and sometimes
simplify our arguments.

\subsubsection{}

Let $\BF$ be an algebraic closure of a finite field and let $X$ be
an algebraic stack of finite type over $\BF$. We denote by $\D(X) =
\D^b_c(X,\bar\Q_l)$ the bounded derived category of constructible
$\bar\Q_l$-sheaves on $X$ (cf. \cite[2.2.18]{BBD82},
\cite[3.3]{Jo93}), where $l$ is a prime number invertible in $\BF$.
An object of $\D(X)$ is also referred to as a complex.

The constant sheaf on $X$ is denoted as $\bar\Q_{l,X}$, or merely
$\bar\Q_l$ when $X$ is clear from the context. We choose an
isomorphism $\bar\Q_{l,\pt}(1)\cong\bar\Q_{l,\pt}$ once and for all
and omit the Tate twist throughout this paper.

\subsubsection{}

Verdier duality and the bifunctors $\RHom,\otimes$ are defined for
$\D(X)$. The Verdier dual of $A\in\D(X)$ is denoted as $DA$.

For a morphism $f: X \to Y$ of algebraic stacks, there is an induced
functor $f^*: \D(Y)\to\D(X)$. When $f$ is representable, there are
also induced functors $f_!,f_*: \D(X)\to\D(Y)$ and $f^!:
\D(Y)\to\D(X)$.

Basic properties (for example, those studies in \cite[Chapter
II,III]{KW01}) of the functors $f^*,f^!,f_*,f_!,D,\RHom,\otimes$ for
schemes have a natural version for algebraic stacks.

\subsubsection{}

For an embedding of locally closed substack $j: S\to X$, we use
$A|_S$ to denote $j^*A$ for a complex $A\in\D(X)$.

\subsubsection{}

Let $\M(X)$ denote the full subcategory of $\D(X)$ consisting of the
perverse sheaves. Let $\pH^n: \D(X)\to\M(X)$ denote the $n$-th
cohomological functor associated to the perverse t-structure.

A complex $C\in\D(X)$ is said to be {\em semisimple} if
$C\cong\oplus_n\pH^n(C)[-n]$ and if $\pH^n(C)\in\M(X)$ is semisimple
for all $n$.

\subsubsection{Decomposition theorem}

Let $f: X\to Y$ be a representable proper morphism of algebraic
stacks with $X$ smooth. Then $f_!\bar\Q_{l,X}\in \D(Y)$ is a
semisimple complex. Cf. \cite[4.13]{Jo93}.

\subsubsection{}\label{sec:pre:perverse:decom}

More generally, let $f: X\to Y$ be a representable morphism of
algebraic stacks. Assume there is a stratification $X=\sqcup_nU_n$
by locally closed substacks such that restricting to each stratum
$U_n$ the morphism $f$ can be factored as
$$U_n \xrightarrow{f'_n} X_n \xrightarrow{f_n} Y$$
where $f'_n$ is a vector bundle of fiber dimension $d_n$ and $f_n$
is a representable proper morphism with $X_n$ smooth. Then
$f_!\bar\Q_{l,X}\in \D(Y)$ is a semisimple complex. Moreover,
$$f_!\bar\Q_{l,X} \cong \oplus_n f_{n!}\bar\Q_{l,X_n}[-2d_n].$$

The proof is essentially the same as \cite[3.7]{Lu85}, which argues
that the complexes $(f|_{U_n})_!\bar\Q_{l,U_n}$ are induced from
pure complexes of the same weight.

\subsection{Fourier-Deligne transform}\label{sec:pre:fourier}

We refer the readers to \cite{KW01} as a general reference for this
subsection.

Let $\BF_p$ be a finite field with $p$ elements, where $p$ is the
characteristic of $\BF$. We fix a nontrivial character $\chi:
\BF_p\to\bar\Q_l^\times$. The Artin-Schreier covering $\BF\to\BF$
given by $x\to x^p-x$ has $\BF_p$ as a group of covering
transformations. Hence the character $\chi$ gives rise to a
$\bar\Q_l$-local system $\CL$ of rank one on the affine line $\BF$.

\subsubsection{}

Let $\pi: E\to X$, $\pi': E'\to X$ be two vector bundles of constant
fiber dimension $d$ over an algebraic variety $X$. Assume we are
given a bilinear map $T: E\times_XE'\to\BF$ which defines a duality
between the two vector bundles. We have a diagram
$E\xleftarrow{s}E\times_XE'\xrightarrow{t}E'$ where $s,t$ are the
obvious projections.

The {\em Fourier-Deligne transform} is defined to be the functor
\begin{eqnarray*}
  && \Phi_{E,E'}: \D(E) \to \D(E'), \quad
  A \mapsto t_!(s^*A\otimes T^*\CL)[d].
\end{eqnarray*}
Interchanging the roles of $E,E'$ we have another Fourier-Deligne
transform
\begin{eqnarray*}
  && \Phi_{E',E}: \D(E') \to \D(E), \quad
  B \mapsto s_!(t^*B\otimes T^*\CL)[d].
\end{eqnarray*}

\subsubsection{Fourier inversion formula}\label{sec:pre:fourier:iso}

There is an isomorphism of functors
$$\Phi_{E',E}\Phi_{E,E'} \cong \sigma^*$$
where $\sigma$ is the multiplication of $-1$ on each fiber of $E$.

It is known that $\Phi_{E,E'}$ is perverse t-exact. In particular,
it restricts to an equivalence of categories $\M(E)\to\M(E')$.

\subsubsection{}\label{sec:pre:fourier:dual}

Fourier-Deligne transform almost commutes with Verdier duality
$$D\Phi_{E,E'}D \cong \Phi_{E,E'}\sigma^*.$$
In fact, by \ref{sec:pre:fourier:iso} the right hand side is inverse
to $\Phi_{E',E}$ hence is right adjoint to $\Phi_{E',E}$. On the
other hand, the left hand side is also right adjoint to
$\Phi_{E',E}$ as demonstrated by the natural isomorphisms
\begin{equation*}
\begin{split}
  &\, \Hom_{\D(E')}(A,D\Phi_{E,E'}DB) \\
  = &\, \Hom_{\D(E')}(A,D[d]t_!(s^*DB\otimes T^*\CL)) \\
  = &\, \Hom_{\D(E')}(A,[-d]t_*\RHom(T^*\CL,s^!B)) \\
  = &\, \Hom_{\D(E)}(s_!(t^*A\otimes T^*\CL)[d],B) \\
  = &\, \Hom_{\D(E)}(\Phi_{E',E}A,B).
\end{split}
\end{equation*}

\subsubsection{}

Let $G$ be an algebraic group. Assume the morphisms
$\pi,\pi',T,\sigma$ are $G$-equivariant (letting $G$ acts trivially
on $\BF$). One may equally well define a Fourier-Deligne transform
$$\Phi_{E',E}: \D([E'/G]) \to \D([E/G])$$
by using the induced morphisms
$[E/G]\xleftarrow{s}[E\times_XE'/G]\xrightarrow{t}[E'/G]$ and $T:
[E\times_XE'/G]\to\BF$. The claims from \ref{sec:pre:fourier:iso}
and \ref{sec:pre:fourier:dual} remain true for this Fourier-Deligne
transform, whose proof is the same as the usual version.

\subsection{Localization in triangulated categories}\label{sec:pre:loc}

A typical example of localization appears in the construction of a
derived category, via which all acyclic complexes are forced to zero
whence quasi-isomorphisms become isomorphisms.

In some sense, the philosophy of localization is that one kills
redundant objects off an ambient category so as to get an
appropriate one for usage.

In the present paper, localization is adopted in the same way as in
micro-local sheaf theory (cf. \cite[6.1]{KSc90}). Objects being
killed are those ``supported outside a given place''. We will
explain this point later in \ref{sec:cat_fun:micro}.

\subsubsection{}

Let $\D$ be a triangulated category. A {\em thick subcategory} of
$\D$ is a full triangulated subcategory $\CN$ such that if a
morphism $A\to B$ in $\D$ factors through an object in $\CN$ and can
be embedded into an exact triangle $A\to B\to C\to A[1]$ with
$C\in\CN$ then $A,B\in\CN$. (In particular, if two objects in an
exact triangle are contained in $\CN$, so is the third one.)

Given a thick subcategory $\CN$, the class of morphisms
\begin{equation*}
\begin{split}
  S := \{s: A\to B \mid & \; \text{$s$ can be embedded into an exact triangle} \\
  & \text{$A\xrightarrow{s}B\to C\to A[1]$ with $C\in\CN$} \}
\end{split}
\end{equation*}
is localizing, hence gives rise to a triangulated category
$\D[S^{-1}]$ by localization (cf. \cite[4.1, 5.1.10]{GM94}). We
denote this category by $\D/\CN$.

By definition, the objects of $\D/\CN$ are the same as $\D$. A
morphism $A\to B$ in $\D/\CN$ is an equivalence class of {\em
roofs}, i.e. an equivalence class of pairs of morphisms
$(T\xrightarrow{s}A,\,T\xrightarrow{t}B)$ in $\D$ with $s\in S$. Two
roofs $(s_i,t_i)$, $i=1,2$ are {\em equivalent}, if there is a third
roof $(s,t)$ forming into a commutative diagram
\begin{equation*}
\xymatrixrowsep{1.5pc}\xymatrixcolsep{.5pc}\xymatrix{
  & & T \ar[ld]_{s} \ar[rd]^{t} \\
  & T_1 \ar[ld]_{s_1} \ar[rrrd]^{t_1} & & T_2 \ar[llld]_{s_2} \ar[rd]^{t_2} \\
  A & & & & B
}
\end{equation*}

\begin{rem}
(1) It can be shown that an object in $\D/\CN$ is isomorphic to zero
if and only if it is an object from $\CN$.

(2) For localization purpose it suffices to assume $\CN$ is a null
system (cf. \cite[1.6]{KSc90}). But this does not provide new
constructions. In fact, a null system $\CN$ can always be completed
to a thick subcategory $\TN$ (the one generated by $\CN$) so that
the localized categories $\D/\CN$, $\D/\TN$ are identical.
\end{rem}

\subsubsection{}\label{sec:pre:loc:thick}

Let $\varphi: \D_1\to\D_2$ be a functor of triangulated categories
and let $\CN_2$ be a thick subcategory of $\D_2$. Then the full
subcategory $\CN_1$ of $\D_1$ whose objects are those sent by
$\varphi$ into $\CN_2$ is a thick subcategory.

\subsubsection{}\label{sec:pre:loc:functor}

Let $\CN_i$ be a thick subcategory of $\D_i$, $i=1,2$ and let
$\varphi:\D_1\to\D_2$ be a functor of triangulated categories such
that $\varphi(\CN_1)\subset\CN_2$. Then the following defines a
functor of triangulated categories
\begin{eqnarray*}
  \tilde\varphi : \D_1/\CN_1 & \to & \D_2/\CN_2 \\
  A & \mapsto & \varphi(A) \\
  ~ [s,t] & \mapsto & [\varphi(s),\varphi(t)].
\end{eqnarray*}

\subsubsection{}\label{sec:pre:loc:transform}

Assume $\psi:\D_1\to\D_2$ satisfies the same assumptions as
$\varphi$. Then a natural transformation $\alpha_A: \varphi A \to
\psi A$, $A\in\D_1$ induces a natural transformation of the induced
functors $[\Id_A,\alpha_A]: \tilde\varphi A \to \tilde\psi A$,
$A\in\D_1/\CN_1$.

It follows that functor isomorphisms and functor adjunctions are
preserved under localization.

\subsubsection{}\label{sec:pre:loc:finite}

Let $(\D^{\ge0},\D^{\le0})$ be a t-structure on $\D$ with the core
$\M=\D^{\ge0}\cap\D^{\le0}$ and cohomological functors $H^n:
\D\to\M$. We assume
\begin{enumerate}
\item
the given t-structure is bounded, i.e. for every object $A\in\D$,
$H^n(A)\not\cong0$ for only finitely many $n$ and, moreover,
$A\cong0$ if and only if $H^n(A)\cong0$ for all $n$;
\item
every object in $\M$ has finite length.
\end{enumerate}

Let $\CN$ be a full subcategory of $\D$ such that
\begin{enumerate}
\setcounter{enumi}{2}
\item
$A\in\CN$ if and only if $H^n(A)\in\CN$ for all $n$;
\item
$\M\cap\CN$ is a Serre subcategory of $\M$, i.e. $\M\cap\CN$ is
stable under extensions and subquotients.
\end{enumerate}
It is easy to check that $\CN$ is the thick subcategory of $\D$
generated by the simple objects from $\M\cap\CN$.

If two thick subcategories $\CN_1,\CN_2$ of $\D$ satisfy (3)(4), so
does the thick subcategory $\TN$ generated by $\CN_1,\CN_2$. In
fact, $\TN$ is generated by the simple objects from $\M\cap\CN_1$
and $\M\cap\CN_2$.

\subsubsection{}\label{sec:pre:loc:t}

With the above assumptions, let $\D/\CN^{\ge0}$ (resp.
$\D/\CN^{\le0}$, $\M/\CN$) be the full subcategory of $\D/\CN$
consisting of the objects isomorphic to those from $\D^{\ge0}$
(resp. $\D^{\le0}$, $\M$). Some labors on homological algebras will
convince the reader that ($\D/\CN^{\ge0},\D/\CN^{\le0})$ defines a
t-structure on $\D/\CN$ with the core $\M/\CN$.

Moreover, an object $A\in\M$ is simple in $\M/\CN$ if and only if
precisely one simple component of its Jordan-H\"older decomposition
in $\M$ is not contained in $\CN$.

\subsubsection{Example}\label{sec:pre:loc:exam}

Let $X$ be an algebraic stack and let $j: U\to X$ be the inclusion
of an open substack. Let $\CN$ be the full subcategory of $\D(X)$
consisting of the objects whose supports are disjoint from $U$.

Then the triangulated category $\D(X)$ endowed with the perverse
t-structure and the full subcategory $\CN$ satisfy the assumptions
of \ref{sec:pre:loc:finite}. (Condition (2) is stated in
\cite[4.3]{Jo93}; the others are clear from definition.)

Moreover, the functor $j^*: \D(X)\to\D(U)$ and the functor $j_!:
\D(U)\to\D(X)$ induce (in the sense of \ref{sec:pre:loc:functor}) an
equivalence of triangulated categories
$$\D(X)/\CN \sim \D(U).$$

\section{Categories and functors from quiver varieties}\label{sec:cat_fun}

Let the quiver $(I,H)$, the generalized Cartan matrix $(a_{ij})$ and
the quantum group $U$ be defined in \ref{sec:pre:u}.

\subsection{The category $\FD$}

\subsubsection{}

First, we enlarge $(I,H)$ to a quiver $(\hI,\hH)$ by appending to
the underlying graph for each vertex $i\in I$ a new vertex $\hi$ and
an edge joining $i,\hi$.
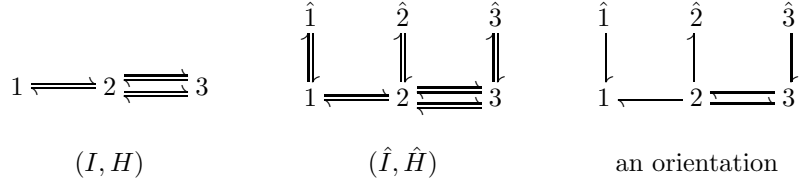
\begin{equation*}
\begin{array}{ccc}
\xymatrixrowsep{1.5pc}\xymatrix{
  \\
  1 \ar@<-.2ex>@_{<-}[r] \ar@<.2ex>@^{->}[r]
  & 2 \ar@<-.9ex>@_{<-}[r] \ar@<-.5ex>@^{->}[r] \ar@<.5ex>@_{<-}[r] \ar@<.9ex>@^{->}[r]
  & 3
} \quad & \quad \xymatrixrowsep{1.5pc}\xymatrix{
  \hat1 \ar@<-.2ex>@_{<-}[d] \ar@<.2ex>@^{->}[d]
  & \hat2 \ar@<-.2ex>@_{<-}[d] \ar@<.2ex>@^{->}[d]
  & \hat3 \ar@<-.2ex>@_{<-}[d] \ar@<.2ex>@^{->}[d] \\
  1 \ar@<-.2ex>@_{<-}[r] \ar@<.2ex>@^{->}[r]
  & 2 \ar@<-.9ex>@_{<-}[r] \ar@<-.5ex>@^{->}[r] \ar@<.5ex>@_{<-}[r] \ar@<.9ex>@^{->}[r]
  & 3
} \quad & \quad \xymatrixrowsep{1.5pc}\xymatrix{
  \hat1 \ar@<.2ex>@^{->}[d]
  & \hat2 \ar@<-.2ex>@_{<-}[d]
  & \hat3 \ar@<.2ex>@^{->}[d] \\
  1 \ar@<-.2ex>@_{<-}[r]
  & 2 \ar@<-.5ex>@^{->}[r] \ar@<.5ex>@_{<-}[r]
  & 3
} \\ \\
(I,H) \quad & \quad (\hI,\hH) \quad & \quad \text{an orientation}
\end{array}
\end{equation*}

\subsubsection{Notations}

Let $\bar{}: \hH\to\hH$ be the involution defined by reversing the
orientation. An {\em orientation} of $(\hI,\hH)$ is a subset
$\Omega\subset\hH$ such that $\Omega\cup\bar\Omega=\hH$ and
$\Omega\cap\bar\Omega=\emptyset$. (In particular, $\Omega$ has the
half number of elements as $\hH$.)

For an arrow $h\in\hH$, we denote by $h',h''$ the source vertex and
target vertex of $h$, respectively. A vertex $i\in\hI$ is called a
{\em source} (resp. {\em sink}) of $\Omega\subset\hH$, if every
arrow $h\in\Omega$ incident to $i$ is outgoing (resp. incoming).

\subsubsection{Quiver varieties}

Assume $\nu=\sum_{i\in\hI}{\nu_i}\cdot i\in\N[\hI]$. We set
$V_i=\BF^{\nu_i}$, $i\in\hI$.

For each subset $\Omega\subset\hH$ we define a variety
$$E_{\Omega} = E_{\nu,\Omega} := \bigoplus_{h\in\Omega} \Hom(V_{h'},V_{h''}).$$
Let the algebraic group
$$\Gn := \prod_{i\in I}GL(V_i)$$
act on $E_{\Omega}$ from the right such that $(x\cdot
g)_h=g^{-1}_{h''}x_hg_{h'}$ (we define $g_{\hi}=\Id$, $i\in I$).

\subsubsection{Fourier-Deligne transform}

Given two orientations $\Omega,\Omega'\subset\hH$, we regard
$E_{\Omega}$ and $E_{\Omega'}$ as vector bundles over
$E_{\Omega\cap\Omega'}$ in the obvious way. The morphism
$$T : E_{\Omega\cup\Omega'}=E_{\Omega} \times_{E_{\Omega\cap\Omega'}} E_{\Omega'} \to \BF, \quad
  x \mapsto \sum_{h\in\Omega\setminus\Omega'}\Tr(x_{\bar h}x_h)
$$
is clearly $\Gn$-equivariant (letting $\Gn$ acts trivially on $\BF$)
and defines a duality of the vector bundles $E_{\Omega}$,
$E_{\Omega'}$, hence gives rise to a Fourier-Deligne transform
$$\Phi_{\Omega,\Omega'}: \D([E_{\Omega}/\Gn]) \to \D([E_{\Omega'}/\Gn]).$$

\begin{prop}\label{prop:d:fourier}
There is an isomorphism of functors
$\Phi_{\Omega_2,\Omega_3}\Phi_{\Omega_1,\Omega_2}\cong\Phi_{\Omega_1,\Omega_3}$
for orientations $\Omega_1,\Omega_2,\Omega_3\subset\hH$.
\end{prop}

Assume $\Omega_a=(\Omega\setminus\Omega'_a)\sqcup\bar\Omega'_a$,
$a=1,2,3$, where $\Omega$ is an orientation and
$\Omega'_1,\Omega'_2,\Omega'_3\subset\Omega$ are disjoint subsets.
It is easy to see
$\Phi_{\Omega_a,\Omega_b}\cong\Phi_{\Omega,\Omega_b}\Phi_{\Omega_a,\Omega}$
for $a\ne b$.

Notice that the morphism $\sigma:
[E_{\Omega}/\Gn]\to[E_{\Omega}/\Gn]$ induced by the multiplication
of $-1$ on a $\Hom$ summand of $E_\Omega$ is isomorphic to the
identity. From \ref{sec:pre:fourier:iso} it follows that
$\Phi_{\Omega_2,\Omega}\Phi_{\Omega,\Omega_2}\cong\Id$.

Therefore, $\Phi_{\Omega_2,\Omega_3}\Phi_{\Omega_1,\Omega_2} \cong
\Phi_{\Omega,\Omega_3}\Phi_{\Omega_2,\Omega}\Phi_{\Omega,\Omega_2}\Phi_{\Omega_1,\Omega}
\cong \Phi_{\Omega_1,\Omega_3}$.

\subsubsection{}\label{sec:cat_fun:cat:tilde}

We denote by $x(i)$ the restriction of $x\in E_\Omega$ to the direct
summand
$$\bigoplus_{h\in\Omega:\;h'=i}\Hom(V_{h'},V_{h''}) = \Hom(V_i,V(i))$$
where
$$V(i):= \bigoplus_{h\in\Omega:\;h'=i}V_{h''}.$$
Define an open subvariety of $E_\Omega$ for each $i\in I$
$$\tE\Oi=\tE\nOi := \{ x\in E_\Omega \mid \Ker x(i)=0 \}.$$

When $i$ is a source of $\Omega$, letting $\dot\Omega$ denote the
subset of $\Omega$ consisting of the arrows not incident to $i$, we
have
$$\tE\Oi/GL(V_i) = \{ (\dx,V) \mid \dx\in E_{\dot\Omega}, \;
  \text{$V\subset V(i)$ is a subspace of dimension $\nu_i$} \}.
$$

\subsubsection{Localization}

We choose for each vertex $i\in I$ an orientation
$\Omega_i\subset\hH$ which has $i$ as a source. Denote by $\CN_i$
the thick subcategory of $\D([E_{\Omega_i}/\Gn])$ consisting of the
complexes whose supports are disjoint from the open substack
$[\tE_{\Omega_i,i}/\Gn]$.

Then for an arbitrary orientation $\Omega\subset\hH$, we denote by
$\CN\nOi$ the thick subcategory of $\D([E_{\Omega}/\Gn])$ consisting
of the complexes sent by $\Phi_{\Omega,\Omega_i}$ into $\CN_i$ (cf.
\ref{sec:pre:loc:thick}). Note that $\CN\nOi$ is independent of the
choice of $\Omega_i$.

Define $\CN\nO$ to be the thick subcategory of
$\D([E_{\Omega}/\Gn])$ generated by $\CN\nOi$, $i\in I$.

\subsubsection{}

By \propref{prop:d:fourier}, the categories $\D([E_{\Omega}/\Gn])$
for various orientation $\Omega$ are related via Fourier-Deligne
transforms, thus determine a unique triangulated category up to
equivalence which we denote as $\D([E_\nu/\Gn])$. Accordingly,
$\CN\nO$ for various $\Omega$ give rise to a thick subcategory
$\CN_\nu$ of $\D([E_\nu/\Gn])$.

Now we are in position to define our main category
$$\FD := \bigoplus_{\nu} \FD_\nu$$
where
$$\FD_\nu:=\D([E_\nu/\Gn])/\CN_\nu.$$

Since Fourier-Deligne transforms are perverse t-exact, perverse
t-structure is well defined on $\D([E_\nu/\Gn])$, as well as on
$\FD$ according to \ref{sec:pre:loc:finite}-\ref{sec:pre:loc:exam}.

According to \ref{sec:pre:fourier:dual} and the isomorphism
$\sigma\cong\Id$ claimed in the proof of \propref{prop:d:fourier},
Verdier duality is well defined on $\D([E_\nu/\Gn])$. It is clear
that $D\CN_\nu\subset\CN_\nu$, hence Verdier duality is well defined
on $\FD$, too.

\subsection{A micro-local point of view}\label{sec:cat_fun:micro}

The category $\D([E_{\Omega}/\Gn])/\CN\nO$ actually can be well
understood from the viewpoint of micro-local sheaf theory. For
simplicity, we assume $(I,H)$ is of finite type.

First, we review in brief Nakajima's quiver variety \cite{Na94,Na98}
associated to $(I,H)$ and $\nu\in\N[\hI]$.

Choose an orientation $\Omega\subset\hH$ and identify $E_{\hH}$ with
$T^*E_\Omega$. Then the full quiver variety $E_{\hH}$ is canonically
a symplectic vector space, on which $\Gn$ acts by symplectic
isomorphisms. So, we have a momentum map $\mu: E_{\hH} \to
Lie(\Gn)^*$ which is assumed vanishing at the origin of $E_{\hH}$.

One can show the algebraic group $\Gn$ acts freely on the variety
$$\mu^{-1}(0)^s := \{ x\in \mu^{-1}(0) \mid \Ker x(i)=0, \; i\in I \},$$
and Nakajima's quiver variety is define to be the quotient
$$\FM_\nu:=\mu^{-1}(0)^s/\Gn.$$
Actually $\mu^{-1}(0)$ and hence $\FM_\nu$ are independent of
$\Omega$.

Now, observe that $\mu^{-1}(0)$ is precisely the union of the
conormal varieties to the $\Gn$-orbits of $E_\Omega$. Intuitively we
may take the quotient stack $[\mu^{-1}(0)/\Gn]$ as the cotangent
bundle to $[E_\Omega/\Gn]$. Also observe that $\FM_\nu$ is an open
substack of $T^*[E_\Omega/\Gn]:=[\mu^{-1}(0)/\Gn]$. (The appropriate
``cotangent bundle'' should be a cotangent complex \cite{Il71} and
there does be a meaningful open immersion of $\FM_\nu$ into the
cotangent complex to $[E_\Omega/\Gn]$ one can work out.)

Therefore, the category $\D([E_{\Omega}/\Gn])/\CN\nO$ is nothing but
the localization of $\D([E_{\Omega}/\Gn])$ by the complexes whose
``micro-supports (i.e. characteristic varieties)'' are disjoint from
$\FM_\nu$.

In fact, similar localization has already appeared in the
construction of stacks of micro-local perverse sheaves
\cite{GMV05,Wa04}. After the terminology therein, we refer to the
perverse sheaves in $\FD_\nu$ provisionally (to the extent maybe
inappropriately) as micro-local perverse sheaves on the quiver
variety $\FM_\nu$.

\subsection{The functors $\FK_i,\FE^{(n)}_i,\FF^{(n)}_i$}\label{sec:cat_fun:fun}

\subsubsection{}\label{sec:cat_fun:fun:kef}

Assume $\nu'=\nu+n\cdot i$ for a vertex $i\in I$ and an integer
$n\ge1$. Let $\Gnn := \Gn\times \Gnp$ act on the variety
$$F_{\nu\nu'} := \{ y\in\bigoplus_{j\in\hI}\Hom(V_j,V'_j) \mid \Ker y_j=0, \; y_{\hj}=\Id, \; j\in I \}$$
from the right such that $(y\cdot(g,g'))_j=g'^{-1}_jy_jg_j$.

We associate to every $\Omega\subset\hH$ a variety
$$Z_\Omega := \{ (x,y) \in E\npO\times F_{\nu\nu'} \mid
  \Img x_hy_{h'}\subset\Img y_{h''} \}
$$
and a couple of morphisms
\begin{equation*}
\begin{array}{ll}
  p: [Z_\Omega/\Gnn] \to [E\nO/\Gn], & p(x,y)_h = y^{-1}_{h''}x_hy_{h'}, \\
  p': [Z_\Omega/\Gnn] \to [E\npO/\Gnp], & p'(x,y)_h = x_h. \\
\end{array}
\end{equation*}

Define the following functors for an orientation $\Omega\subset\hH$
\begin{equation*}
\begin{array}{l@{\;}l}
  \K\Oi:=\Id[\nu_i-\bar\nu_i] & : \D([E\nO/\Gn])\to\D([E\nO/\Gn]), \\
  \F^{(n)}\Oi:=p'_!p^*[n\bar\nu'_i(\Omega)] & : \D([E\nO/\Gn])\to\D([E\npO/\Gnp]). \\
\end{array}
\end{equation*}
where
\begin{equation*}
\begin{split}
  & \bar\nu_i := \sum_{h\in\hH: \; h'=i}\nu_{h''}-\nu_i, \\
  & \bar\nu_i(\Omega) := \sum_{h\in\Omega: \; h'=i} \nu_{h''}-\nu_i. \\
\end{split}
\end{equation*}

When $i$ is a source of $\Omega\subset\hH$, we have a subvariety of
$Z_\Omega$
$$\tZ\Oi := \{ (x,y) \in \tE\npOi\times F_{\nu\nu'} \mid
  \Img x_hy_{h'}\subset\Img y_{h''} \}.
$$
Let $\tp_i,\tp'_i$ be the restrictions of $p,p'$ to $[\tZ\Oi/\Gnn]$,
respectively, and define
$$\E^{(n)}\Oi:=\tp_{i!}\tp'^*_i[n\nu_i] : \D([E\npO/\Gnp])\to\D([E\nO/\Gn]).$$

\begin{rem}
(1) A $G$-equivariant morphism $X\to Y$ of algebraic varieties
automatically induces a representable morphism $[X/G]\to[Y/G]$. If,
in addition, $G$ acts trivially on $Y$, there is also an induced
(maybe non-representable) morphism $[X/G]\to Y$. As we have done for
$p,p'$ above, we usually describe induced morphisms of quotient
stacks by their liftings.

(2) The morphisms $p',\tp_i,\tp'_i$ are representable and smooth,
while $p$ may be even not representable. But the functor $p^*$ is
well defined.

(3) It is not surprising to see the definitions of $\E, \F$ are less
symmetric, for we are in purpose to treat highest weight
representations of $U$.
\end{rem}

\begin{prop}\label{prop:kef:fourier}
There are functor isomorphisms for orientations
$\Omega_1,\Omega_2\subset\hH$ (both having $i$ as a source for the
isomorphism of $\E$)
\begin{equation*}
\begin{split}
  \Phi_{\Omega_1,\Omega_2}\K_{\Omega_1,i}& \cong\K_{\Omega_2,i}\Phi_{\Omega_1,\Omega_2}, \\
  \Phi_{\Omega_1,\Omega_2}\E^{(n)}_{\Omega_1,i}& \cong\E^{(n)}_{\Omega_2,i}\Phi_{\Omega_1,\Omega_2}, \\
  \Phi_{\Omega_1,\Omega_2}\F^{(n)}_{\Omega_1,i}& \cong\F^{(n)}_{\Omega_2,i}\Phi_{\Omega_1,\Omega_2}.
\end{split}
\end{equation*}
\end{prop}

The isomorphism of $\K$ is obvious. The proof of $\E$ is similar as
$\F$. Below we prove the isomorphism of $\F$.

We may assume $\Omega_1,\Omega_2$ differ by a single element. That
is, $\Omega_1\setminus\Omega_2=\{h\}$. We form the following
commutative diagram.
\begin{equation*}
\xymatrixcolsep{1.5pc}\xymatrix{
  [E_{\nu,\Omega_2}/\Gn]
  & [Z_{\Omega_2}/\Gnn] \ar[l]_--p \ar[r]^--{p'}
  & [E_{\nu',\Omega_2}/\Gnp]\\
  [E_{\nu,\Omega_1\cup\Omega_2}/\Gn] \ar[u]_{t} \ar [d]^{s}
  & [Z_{\Omega_1\cup\Omega_2}/\Gnn] \ar[l]_--p \ar[r]^--{p'} \ar[u]_{t_Z} \ar [d]^{s_Z}
  & [E_{\nu',\Omega_1\cup\Omega_2}/\Gnp] \ar[u]_{t'} \ar [d]^{s'} \\
  [E_{\nu,\Omega_1}/\Gn]
  & [Z_{\Omega_1}/\Gnn] \ar[l]_--p \ar[r]^--{p'}
  & [E_{\nu',\Omega_1}/\Gnp]
}
\end{equation*}
Let $T_Z$ be the morphism
$$T_Z: [Z_{\Omega_1\cup\Omega_2}/\Gnn] \to \BF, \quad
  (x,y) \mapsto \Tr(x_{\bar h}x_h).
$$

When $h$ is not incident to $i$, all the squares in the above
diagram are cartesian. Otherwise, we assume $h''=i$. Then the
top-left and the bottom-right squares in the above diagram are
cartesian. In either case, there are natural isomorphisms for
$A\in\D([E_{\nu,\Omega_1}/\Gn])$
\begin{equation*}
\begin{split}
  p^*\Phi_{\Omega_1,\Omega_2}A[-\nu_{h'}\nu_{h''}]
  & = p^*t_!(s^*A\otimes T^*\CL)
  \cong t_{Z!}p^*(s^*A\otimes T^*\CL) \\
  & \cong t_{Z!}(p^*s^*A\otimes T_Z^*\CL)
  = t_{Z!}(s_Z^*p^*A\otimes T_Z^*\CL).
\end{split}
\end{equation*}
Similarly, there are natural isomorphisms for
$B\in\D([Z_{\Omega_1}/\Gnn])$
$$\Phi_{\Omega_1,\Omega_2}p'_!B[-\nu'_{h'}\nu'_{h''}] \cong p'_!t_{Z!}(s_Z^*B\otimes T_Z^*\CL).$$
It follows that
\begin{equation*}
  p'_!p^*\Phi_{\Omega_1,\Omega_2}[-\nu_{h'}\nu_{h''}]
  \cong \Phi_{\Omega_1,\Omega_2}p'_!p^*[-\nu'_{h'}\nu'_{h''}].
\end{equation*}
Therefore, the isomorphism of $\F$ follows.

\begin{prop}\label{prop:kef:local}
Let $\Omega\subset\hH$ be an orientation (having $i$ as a source for
the statement of $\E$). Then $\K\Oi(\CN\nO)\subset\CN\nO$,
$\E^{(n)}\Oi(\CN\npO)\subset\CN\nO$,
$\F^{(n)}\Oi(\CN\nO)\subset\CN\npO$.
\end{prop}

If $k\in I$ is a source of $\Omega$, we have clearly
$\K\Oi(\CN\nOk)\subset\CN\nOk$,
$\E^{(n)}\Oi(\CN\npOk)\subset\CN\nOk$,
$\F^{(n)}\Oi(\CN\nOk)\subset\CN\npOk$. Hence from the above
proposition the claim for $\K,\F$ follows.

We are left to show $\E^{(n)}\Oi(\CN\npOk)\subset\CN\nO$ where $k\in
I$ is a vertex joined by some edge(s) to $i$. Below we assume $n=1$,
proof of general case postponed to the next subsection.

Let $\Omega'$ be obtained from $\Omega$ by reversing the arrows
joining $k,i$. We may assume $k$ is a source of $\Omega'$. We have
to show \setno
\begin{equation}\label{eqn:kef:local1}
  \tp_{i!}\tp'^*_iA\in\CN\nO
  \;\;\;\text{for}\;\;\; A\in\CN\npOk.
\end{equation}

Let $H_0\subset H$ be the set of the arrows going from $k$ to $i$.
Note that $\Omega\cup\Omega' = \Omega\sqcup H_0 = \Omega'\sqcup\bar
H_0$. Define for every subset $\tilde\Omega\subset\hH$ a variety
$$W_{\tilde\Omega} := \{ (x,x',y,y') \in
  E_{\nu,\tilde\Omega\setminus H_0}\times E_{\nu',\tilde\Omega\cap H_0}\times F_{\nu\nu'}\times \bigoplus_{j\in\hI}\Hom(V'_j,V_j)
  \mid y'_jy_j=\Id_{V_j} \}
$$
and a couple of morphisms
\begin{equation*}
\begin{array}{ll}
  q: [W_{\tilde\Omega}/\Gnn] \to [E_{\nu,\tilde\Omega}/\Gn],
  & \text{$q(x,x',y,y')_h = x_h$ or $y'_{h''}x'_hy_{h'}$}, \\
  q': [W_{\tilde\Omega}/\Gnn] \to [E_{\nu',\tilde\Omega}/\Gnp],
  & \text{$q'(x,x',y,y')_h = y_{h''}x_hy'_{h'}$ or $x'_h$}. \\
\end{array}
\end{equation*}

Assume $A\cong\Phi_{\Omega',\Omega}B$ where
$B\in\CN_{\nu'\Omega',k}$. By the same reason as
\propref{prop:kef:fourier}, we have $\Phi_{\Omega',\Omega}p_!p'^*B
\cong p_!p'^*\Phi_{\Omega',\Omega}B [n_1]$ and
$\Phi_{\Omega',\Omega}q_!q'^*B \cong q_!q'^*\Phi_{\Omega',\Omega}B
[n_2]$ for some integers $n_1,n_2$. Clearly $p_!p'^*B, \; q_!q'^*B
\in \CN_{\nu,\Omega',k}$. Hence \setno
\begin{equation}\label{eqn:kef:local2}
  p_!p'^*A, \; q_!q'^*A \in \CN\nOk.
\end{equation}

We stratify $Z_\Omega,W_\Omega$ as follows
\begin{equation*}
\begin{split}
  Z_1 & := \{ (x,y)\in Z_\Omega \mid p(x,y)\notin\tE\nOi \}, \\
  Z_2 & := Z_3\setminus Z_1, \\
  Z_3 & := Z_\Omega\setminus Z_4, \\
  Z_4 & := \{ (x,y)\in Z_\Omega \mid x\in\tE\npOi \}, \\
  W_1 & := \{ (x,,y,y')\in W_\Omega \mid q(x,,y,y')\notin\tE\nOi \}, \\
  W_2 & := W_\Omega\setminus W_1, \\
\end{split}
\end{equation*}
and denote by $z_a,w_a$ the inclusions of locally closed substacks
\begin{equation*}
\begin{split}
  z_a & : [Z_a/\Gnn] \to [Z_\Omega/\Gnn], \\
  w_a & : [W_a/\Gnn]\to [W_\Omega/\Gnn].
\end{split}
\end{equation*}

By definition, \setno
\begin{equation}
  p_!z_{1!}z_1^*p'^*A, \; q_!w_{1!}w_1^*q'^*A \in \CN\nOi.
\end{equation}
Moreover, our assumption $n=1$ implies the morphism
$$[W_2/\Gnn]\to[Z_2/\Gnn],\quad (x,,y,y')\mapsto(yxy',y)$$
is an isomorphism which forms into the following commutative
diagram.
\begin{equation*}
\xymatrixcolsep{1.5pc}\xymatrix{
  [Z_2/\Gnn] \ar[r]^--{p'z_2} \ar[d]_{pz_2} & [E\npO/\Gn] \\
  [E\nO/\Gn] & \ar[l]^--{qw_2} [W_2/\Gnn] \ar[u]_{q'w_2} \ar[lu]_\sim \\
}
\end{equation*}
Thus \setno
\begin{equation}\label{eqn:kef:local4}
  q_!w_{2!}w_2^*q'^*A \cong p_!z_{2!}z_2^*p'^*A.
\end{equation}

From \eqref{eqn:kef:local2}-\eqref{eqn:kef:local4}, it follows that
all the complexes in the following adjunction triangles are sent by
$p_!$ or $q_!$ into $\CN\nO$.
\begin{equation*}
\begin{array}{c@{~\to~}c@{~\to~}c@{~\xrightarrow{[1]}}}
  z_{4!}z_4^*p'^*A & p'^*A & z_{3!}z_3^*p'^*A \\
  z_{2!}z_2^*p'^*A & z_{3!}z_3^*p'^*A & z_{1!}z_1^*p'^*A \\
  w_{2!}w_2^*q'^*A & q'^*A & w_{1!}w_1^*q'^*A \\
\end{array}
\end{equation*}

Finally, observe that $\tp_{i!}\tp'^*_iA \cong p_!z_{4!}z_4^*p'^*A$.
Our claim \eqref{eqn:kef:local1} follows.

\begin{rem}
We have cheated in the above proof: the functors $p_!,q_!$ are not
well defined since the morphisms $p,q$ are not representable.

We can overcome this gap as follows. Choose a connected smooth
principal $GL(V'_i)$-bundle $P$ and replace the functors
$\tp_{i!},p_!,q_!$ by the compositions of the following well defined
ones (defined in the obvious way)
\begin{equation*}
\begin{array}{c@{\;}c@{\;}c@{\;}c@{\;}c}
  \D([\tZ\Oi/\Gnn]) & \xrightarrow{s^*} & \D([\tZ\Oi\times P/\Gnn]) & \xrightarrow{(\tp_is)_!} & \D([E\nO/\Gn]) \\
  \D([Z_\Omega/\Gnn]) & \xrightarrow{s^*} & \D([Z_\Omega\times P/\Gnn]) & \xrightarrow{(ps)_!} & \D([E\nO/\Gn]) \\
  \D([W_\Omega/\Gnn]) & \xrightarrow{s^*} & \D([W_\Omega\times P/\Gnn]) & \xrightarrow{(qs)_!} & \D([E\nO/\Gn]) \\
\end{array}
\end{equation*}
Then the above proof goes through to yield
\begin{equation*}
  \tp_{i!}s_!s^*\tp'^*_iA\in\CN\nO
  \;\;\;\text{for}\;\;\; A\in\CN\npOk.
\end{equation*}

Now, given a perverse sheaf $A\in\CN\npOk$, let $m$ be such that
$\ptau_{\le m}\ptau_{\ge-m}\tp_{i!}\tp'^*_iA \cong
\tp_{i!}\tp'^*_iA$ where $\ptau_{\le m},\ptau_{\ge-m}$ are the
truncation functors associated to the perverse t-structure. Choosing
$P$ such that $H^n(P,\bar\Q_l)=0$ for $1\le n\le2m$, we have
\begin{equation*}
\begin{split}
  & \ptau_{\ge-m}(\tp_{i!}s_!s^*\tp'^*_iA[2\dim P])
  \cong \ptau_{\ge-m}\tp_{i!} \big( s_![2\dim P]\bar\Q_l\otimes\tp'^*_iA \big)
  \cong \tp_{i!}\tp'^*_iA.
\end{split}
\end{equation*}
Since both thick subcategories $\CN\npOk,\CN\nO$ are generated by
simple perverse sheaves therein (cf. \ref{sec:pre:loc:finite},
\ref{sec:pre:loc:exam}), it follows that
$\tp_{i!}\tp'^*_iA\in\CN\nO$ holds for arbitrary $A\in\CN\npOk$.
\end{rem}

\subsubsection{}

By the above propositions, the functors
$\K\Oi,\E^{(n)}\Oi,\F^{(n)}\Oi$ pass to localized categories and
induce the following functors unambiguously
\begin{equation*}
\begin{split}
  \FK_{\nu,i} & : \FD_\nu\to\FD_\nu, \\
  \FE^{(n)}_{\nu,i} & : \FD_{\nu+ni}\to\FD_\nu, \\
  \FF^{(n)}_{\nu,i} & : \FD_\nu\to\FD_{\nu+ni}.
\end{split}
\end{equation*}

Assemble these functors to endofunctors of $\FD$
\begin{eqnarray*}
  && \FK_i := \bigoplus_\nu\FK_{\nu,i}, \\
  && \FE^{(n)}_i := \bigoplus_\nu\FE^{(n)}_{\nu,i}, \\
  && \FF^{(n)}_i := \bigoplus_\nu\FF^{(n)}_{\nu,i}.
\end{eqnarray*}

Below we will drop the superscript $(n)$ off when $n=1$.

\subsection{An alternative approach to $\FE^{(n)}_i,\FF^{(n)}_i$}\label{sec:cat_fun:alt}

Recall that for every vertex $k\in I$ there is an open subvariety
$\tE\nOk\subset E\nO$. Let $j\Ok : [\tE\nOk/\Gn] \to [E\nO/\Gn]$ be
the induced open immersion.

\subsubsection{}

Let $\Omega\subset\hH$ be an orientation having a vertex $k\in I$ as
a source. There are a couple of morphisms obtained from $p,p'$ by
restriction
\begin{equation*}
\begin{array}{ll}
  p_k : [\tZ\Ok/\Gnn] \to [\tE\nOk/\Gn], & p_k(x,y) = y^{-1}xy, \\
  p'_k : [\tZ\Ok/\Gnn] \to [\tE\npOk/\Gnp], & p'_k(x,y) = x. \\
\end{array}
\end{equation*}
We define a functor
\begin{equation*}
  \TF^{(n)}_{\Omega,k,i}:=p'_{k!}p_k^*[n\bar\nu'_i(\Omega)] :
  \D([\tE\nOk/\Gn])\to\D([\tE\npOk/\Gnp]).
\end{equation*}

Assume $k\ne i$. In the notations from \ref{sec:cat_fun:cat:tilde},
we can identify
\begin{equation*}
\begin{split}
  \tZ\Ok/GL(V_k)\times GL(V'_k)
  = & \; \{ (\dx,V,\dy) \mid
  \dx\in E_{\nu',\dot\Omega}, \; \dy\in\dot F_{\nu\nu'}, \; \\
  & \Img \dx_h\dy_{h'}\subset\Img \dy_{h''}, \; V\subset V(k), \; \dim V=\nu_k \}
\end{split}
\end{equation*}
where
$$\dot F_{\nu\nu'} := \{ y\in\bigoplus_{j\in\hI\setminus\{k\}}\Hom(V_j,V'_j) \mid \Ker y_j=0, \; y_{\hj}=\Id, \; j\in I \}$$
so that $p_k(\dx,V,\dy)=(\dy^{-1}\dx\dy,V)$,
$p'_k(\dx,V,\dy)=(\dx,\dy(V))$.

\subsubsection{}\label{sec:cat_fun:alt:p}

Let $\Omega\subset\hH$ be an orientation having $i$ as a source.
Define functors
\begin{equation*}
\begin{array}{l@{\;}l}
  \TE^{(n)}\Oi := p_{i!}p'^*_i[n\nu_i] & : \D([\tE\npOi/\Gnp])\to\D([\tE\nOi/\Gn]), \\
  \TF^{(n)}\Oi := p'_{i!}p^*_i[n\bar\nu'_i] & : \D([\tE\nOi/\Gn])\to\D([\tE\npOi/\Gnp]). \\
\end{array}
\end{equation*}
Note that the functors $\TF^{(n)}_{\Omega,i,i}$, $\TF^{(n)}\Oi$
coincide.

In the notations from \ref{sec:cat_fun:cat:tilde}, we can identify
\begin{equation*}
\begin{split}
  \tZ\Oi/\Gn\times GL(V'_i)
  = & \; \{ (\dx,V\subset V') \mid \dx\in E_{\dot\Omega}, \; V,V'\subset V(i), \; \\
  & \dim V=\nu_i, \; \dim V'=\nu'_i \}
\end{split}
\end{equation*}
so that $p_i(\dx,V,V')=(\dx,V)$, $p'_i(\dx,V,V')=(\dx,V')$.

It follows that $p_i, p'_i$ are fibred in the Grassmannians
$Gr(\bar\nu'_i,\bar\nu_i)$, $Gr(\nu_i,\nu'_i)$ of $\bar\nu'_i,
\nu_i$-dimensional subspaces in a $\bar\nu_i, \nu'_i$-dimensional
$\BF$-vector space, respectively. In particular, the morphisms
$p_i,p'_i$ are proper and smooth of relative dimension
$n\bar\nu'_i$, $n\nu_i$, respectively.

\begin{prop}\label{prop:kef:alt}
We have functor isomorphisms
$$\E^{(n)}\Oi \cong j_{\Omega,i!} \TE^{(n)}\Oi j\Oi^*
  \;\;\;\text{and}\;\;\;
  j\Ok^* \F^{(n)}\Oi \cong \TF^{(n)}_{\Omega,k,i} j\Ok^*.
$$
\end{prop}

The first one follows readily from definition. To see the second
one, we form the following commutative diagram
\begin{equation*}
\xymatrixcolsep{1.5pc}\xymatrix{
  [\tE\nOk/\Gn] \ar[d]_{j\Ok}
  & [\tZ\Ok/\Gnn] \ar[l]_--{p_k} \ar[r]^--{p'_k} \ar[d]^j
  & [\tE\npOk/\Gnp] \ar[d]^{j\Ok} \\
  [E\nO/\Gn]
  & [Z_\Omega/\Gnn] \ar[l]_--p \ar[r]^--{p'}
  & [E\npO/\Gnp] \\
}
\end{equation*}
in which the right square is cartesian. Then
$$j\Ok^*p'_!p^* \cong p'_{k!}j^*p^* = p'_{k!}p_k^*j\Ok^*.$$
Hence
$$j\Ok^*\F^{(n)}\Ok
  = j\Ok^*p'_!p^*[n\bar\nu'_i(\Omega)]
  \cong p'_{k!}p_k^*j\Ok^*[n\bar\nu'_i(\Omega)]
  = \TF^{(n)}_{\Omega,k,i} j\Ok^*.
$$

\begin{prop}\label{prop:kef:EnFn}
We have functor isomorphisms
\begin{align*}
  \TE\Oi^{(n-1)} \TE\Oi & \cong \bigoplus\limits_{0\le m<n} \TE\Oi^{(n)} [n-1-2m], \\
  \TF\Oi^{(n-1)} \TF\Oi & \cong \bigoplus\limits_{0\le m<n} \TF\Oi^{(n)} [n-1-2m].
\end{align*}
\end{prop}

We prove the second isomorphism. Assume $\nu'=\nu+i$,
$\nu''=\nu+ni$. Consider the following commutative diagram which
contains a cartesian square at the bottom-right corner.
$$\xymatrixcolsep{1.5pc}\xymatrix{
  [\tZ^3\Oi/G_{\nu\nu''}] \ar[rr]^{p'_3} \ar[dd]_{p_3} & & [\tE_{\nu'',\Omega,i}/G_{\nu''}] \\
  & [Y/G_{\nu\nu'\nu''}] \ar[lu]_{q_3} \ar[r]^{q_2} \ar[d]^{q_1} & [\tZ^2\Oi/G_{\nu'\nu''}] \ar[u]_{p'_2} \ar[d]^{p_2} \\
  [\tE\nOi/\Gn] & \ar[l]_{p_1} [\tZ^1\Oi/\Gnn] \ar[r]^{p'_1} & [\tE\npOi/\Gnp]
}
$$
where (1) $\tZ^a\Oi,p_a,p'_a$, $a=1,2,3$ are those defining the
functors $\TF\Oi$, $\TF^{(n-1)}\Oi$, $\TF^{(n)}\Oi$, respectively,
(2) $Y := \tZ^1\Oi \times_{\tE\npOi} \tZ^2\Oi$ and
$q_3(\dx,V,V',V'') := (\dx,V,V'')$ under the identification
\begin{equation*}
\begin{split}
  Y/\Gnn\times GL(V''_i)
  = & \; \{ (\dx,V\subset V'\subset V'') \mid
  \dx\in E_{\dot\Omega}, \; V,V',V''\subset V(i), \; \\
  & \dim V=\nu_i, \; \dim V'=\nu'_i, \; \dim V''=\nu''_i \}.
\end{split}
\end{equation*}

We have
\begin{align*}
  \TF\Oi^{(n-1)} \TF\Oi
  = p'_{2!} p_2^* p'_{1!} p_1^*[t]
  \cong p'_{2!} q_{2!} q_1^* p_1^*[t]
  = p'_{3!} q_{3!} q_3^* p_3^*[t]
\end{align*}
where
$$t := \bar\nu'_i+(n-1)\bar\nu''_i = n\bar\nu''_i+n-1.$$

Note that $q_3$ is a $\BP^{n-1}$-bundle. By the decomposition
theorem,
$$q_{3!}\bar\Q_{l,[Y/G_{\nu\nu'\nu''}]} \cong \bigoplus_{0\le m<n}\bar\Q_{l,[\tZ^3\Oi/G_{\nu\nu''}]}[-2m].$$
Hence
$$q_{3!}q_3^* \cong q_{3!}\bar\Q_{l,[Y/G_{\nu\nu'\nu''}]}\otimes- \cong \bigoplus_{0\le m<n}[-2m].$$
It follows that
\begin{equation*}
  \TF\Oi^{(n-1)} \TF\Oi
  \cong \bigoplus\limits_{0\le m<n} p'_{3!} p_3^* [t-2m]
  = \bigoplus\limits_{0\le m<n} \TF\Oi^{(n)} [n-1-2m].
\end{equation*}

\subsubsection{Proof of \propref{prop:kef:local} (continue)}

Since $j^*\Oi j_{\Omega,i!}\cong\Id$, it follows from the above
propositions that
$$\E\Oi^{(n-1)} \E\Oi \cong \bigoplus\limits_{0\le m<n} \E\Oi^{(n)} [n-1-2m].$$
Thus by induction on $n$, the claim
$\E^{(n)}\Oi(\CN\npO)\subset\CN\nO$ is reduced to the case $n=1$,
which has already been proved.

\subsubsection{}

Let $\TN\nOk$ be the thick subcategory of $\D([\tE\nOk/\Gn])$
consisting of those complexes sent by the functor $j_{\Omega,k!}$
into the thick subcategory $\CN\nO$ of $\D([E\nO/\Gn])$ (cf.
\ref{sec:pre:loc:thick}). Then the functors $j_{\Omega,k!}$ and
$j\Ok^*$ induce an equivalence of categories (cf.
\ref{sec:pre:loc:exam})
$$\D([\tE\nOk/\Gn])/\TN\nOk \sim \D([E\nO/\Gn])/\CN\nO.$$
Therefore, $\D([\tE\nOk/\Gn])/\TN\nOk$ provides an alternative
realization of $\FD_\nu$.

Since $j\Ok^*$ is perverse t-exact and commutes with Verdier
duality, the perverse t-structure and Verdier duality on
$\D([E\nO/\Gn])$ induce those on $\FD_\nu$.

Moreover, by \propref{prop:kef:alt},
$\TE^{(n)}\Oi(\TN\npOi)\subset\TN\nOi$,
$\TF^{(n)}_{\Omega,k,i}(\TN\nOk)\subset\TN\npOk$ so that
$\TE^{(n)}\Oi$, $\TF^{(n)}_{\Omega,k,i}$ pass to localized
categories and realize the functors $\FE^{(n)}_{\nu,i}$ and
$\FF^{(n)}_{\nu,i}$, respectively.

\subsection{Functor isomorphisms}\label{sec:cat_fun:iso}

We categorify in this subsection the defining relations of $U$.

\begin{prop}\label{prop:kef:adjoint}
The endofunctors $\FK_i, \FE^{(n)}_i, \FF^{(n)}_i$ of $\FD$ have the
functors
$$\FK_i^{-1}, \quad \FK_i^n\FF^{(n)}_i[-n^2], \quad \FK_i^{-n}\FE^{(n)}_i[-n^2]$$
as left adjoints and have the functors
$$\FK_i^{-1}, \quad \FK_i^{-n}\FF^{(n)}_i[n^2], \quad \FK_i^n\FE^{(n)}_i[n^2]$$
as right adjoints, respectively.
\end{prop}

Recall from \ref{sec:cat_fun:alt:p} that the morphisms $p_i,p'_i$
are proper and smooth of relative dimension $n\bar\nu'_i$, $n\nu_i$,
respectively. Hence $p_{i!}=p_{i*}$ and
$p'^!_i[-n\nu_i]=p'^*_i[n\nu_i]$. We have natural isomorphisms for
$A\in\D([\tE\nOi/\Gn]), B\in\D([\tE\npOi/\Gnp])$
\begin{equation*}
\begin{split}
  & \Hom_{\D([\tE\npOi/\Gnp])} (\TF\Oi^{(n)}A,B) \\
  = & \Hom_{\D([\tE\npOi/\Gnp])} (p'_{i!}p_i^*[n\bar\nu'_i]A,B) \\
  = & \Hom_{\D([\tE\nOi/\Gn])} (A,p_{i*}p'^!_i[-n\bar\nu'_i]B) \\
  = & \Hom_{\D([\tE\nOi/\Gn])} (A,p_{i!}p'^*_i[2n\nu_i-n\bar\nu'_i]B) \\
  = & \Hom_{\D([\tE\nOi/\Gn])} (A,[n\nu_i-n\bar\nu_i]\TE\Oi^{(n)}[n^2]B).
\end{split}
\end{equation*}
It follows that $\FF^{(n)}_i$ is left adjoint to
$\FK_i^n\FE^{(n)}_i[n^2]$ (cf. \ref{sec:pre:loc:transform}).

Similarly for the others.

\begin{thm}\label{thm:kef}
There are isomorphisms of endofunctors of $\FD$.
\begin{enumerate}
\setlength{\itemsep}{1.25ex}
\item
  $\FK_i\FK_j=\FK_j\FK_i$;
\item
  $\FK_i\FE^{(n)}_j = \FE^{(n)}_j\FK_i[-na_{ij}]$;
\item
  $\FK_i\FF^{(n)}_j = \FF^{(n)}_j\FK_i[na_{ij}]$;
\item
  $\FE^{(n-1)}_i \FE_i \cong \bigoplus\limits_{0\le m<n} \FE^{(n)}_i [n-1-2m]$;
\item
  $\FF^{(n-1)}_i \FF_i \cong \bigoplus\limits_{0\le m<n} \FF^{(n)}_i [n-1-2m]$;
\item
  $\FE_i\FF_i \oplus
  \bigoplus\limits_\nu \bigoplus\limits_{0\le m<\nu_i-\bar\nu_i} \Id_{\FD_\nu}[\nu_i-\bar\nu_i-1-2m]$
  \vspace{1.25ex} \\
  $\cong \FF_i\FE_i \oplus
  \bigoplus\limits_\nu \bigoplus\limits_{0\le m<\bar\nu_i-\nu_i} \Id_{\FD_\nu}[\bar\nu_i-\nu_i-1-2m]$;
\item
  $\FE_i\FF_j \cong \FF_j\FE_i, \quad i\ne j$;
\item
  $\bigoplus\limits_{\stackrel{0\le m\le 1-a_{ij}}{\text{$m$ odd}}} \FE^{(m)}_i\FE_j\FE^{(1-a_{ij}-m)}_i
  \cong \bigoplus\limits_{\stackrel{0\le m\le 1-a_{ij}}{\text{$m$ even}}} \FE^{(m)}_i\FE_j\FE^{(1-a_{ij}-m)}_i, \quad i\ne j$;
\item
  $\bigoplus\limits_{\stackrel{0\le m\le 1-a_{ij}}{\text{$m$ odd}}} \FF^{(m)}_i\FF_j\FF^{(1-a_{ij}-m)}_i
  \cong \bigoplus\limits_{\stackrel{0\le m\le 1-a_{ij}}{\text{$m$ even}}} \FF^{(m)}_i\FF_j\FF^{(1-a_{ij}-m)}_i, \quad i\ne j$.
\end{enumerate}
\end{thm}

(1)(2)(3) are obvious. (4)(5) are immediate from
\propref{prop:kef:EnFn}. (8) follows from the above proposition and
(3)(9). The remaining three are proved below.

\subsubsection{Proof of (6)}

Let $\Omega\subset\hH$ be an orientation having $i$ as a source. It
suffices to show there is an isomorphism of endofunctors of
$\D([\tE\nOi/\Gn])$ \setno
\begin{equation}\label{eqn:kef:iso61}
\begin{split}
  & \TE\Oi\TF\Oi \oplus
  \bigoplus_{0\le m<\nu_i-\bar\nu_i} \Id[\nu_i-\bar\nu_i-1-2m] \\
  & \cong \TF\Oi\TE\Oi \oplus
  \bigoplus_{0\le m<\bar\nu_i-\nu_i} \Id[\bar\nu_i-\nu_i-1-2m].
\end{split}
\end{equation}

Assume $\nu^1=\nu+i$, $\nu^2=\nu-i$, $\nu'=\nu$. Consider the
following commutative diagrams, both having a cartesian square at
the bottom-right corner
$$\xymatrixcolsep{1.5pc}\xymatrix{
  [Y/\Gnn] \ar[rr]^{\pi'} \ar[dd]_{\pi} & & [\tE\npOi/\Gnp] \\
  & [Y^1/G_{\nu\nu^1\nu'}] \ar[lu]_{r_1} \ar[r] \ar[d] & [\tZ^1\Oi/G_{\nu'\nu^1}] \ar[u]_{p_1} \ar[d]^{p'_1} \\
  [\tE\nOi/\Gn] & \ar[l]_{p_1} [\tZ^1\Oi/G_{\nu\nu^1}] \ar[r]^{p'_1} & [\tE_{\nu^1,\Omega,i}/G_{\nu^1}]
}
$$
$$\xymatrixcolsep{1.5pc}\xymatrix{
  [Y/\Gnn] \ar[rr]^{\pi'} \ar[dd]_{\pi} & & [\tE\npOi/\Gnp] \\
  & [Y^2/G_{\nu\nu^2\nu'}] \ar[lu]_{r_2} \ar[r] \ar[d] & [\tZ^2\Oi/G_{\nu^2\nu'}] \ar[u]_{p'_2} \ar[d]^{p_2} \\
  [\tE\nOi/\Gn] & \ar[l]_{p'_2} [\tZ^2\Oi/G_{\nu^2\nu}] \ar[r]^{p_2} & [\tE_{\nu^2,\Omega,i}/G_{\nu^2}]
}
$$
where (1) $\tZ^a\Oi,p_a,p'_a$, $a=1,2$ are those defining the
functors $\TE\Oi$, $\TF\Oi$, (2) $Y$ is such that
$$Y/\Gnp\times GL(V_i)
  = \{ (\dx,V,V') \mid \dx\in E_{\dot\Omega}, \; V,V'\subset V(i), \; \\
  \dim V=\dim V'=\nu_i \}
$$
and $\pi(\dx,V,V'):=(\dx,V)$, $\pi'(\dx,V,V'):=(\dx,V')$, (3) $Y^a
:= \tZ^a\Oi \times_{\tE_{\nu^a,\Omega,i}} \tZ^a\Oi$ and
$r_a(\dx,V,V^a,V') := (\dx,V,V')$ under the identification
\begin{equation*}
\begin{split}
  Y^1/G_{\nu^1\nu'}\times GL(V_i)
  = & \; \{ (\dx,V\subset V^1\supset V') \mid
  \dx\in E_{\dot\Omega}, \; V,V^1,V'\subset V(i), \; \\
  & \dim V=\dim V'=\nu_i, \; \dim V^1=\nu_i+1 \}, \\
  Y^2/G_{\nu^2\nu'}\times GL(V_i)
  = & \; \{ (\dx,V\supset V^2\subset V') \mid
  \dx\in E_{\dot\Omega}, \; V,V^2,V'\subset V(i), \; \\
  & \dim V=\dim V'=\nu_i, \; \dim V^2=\nu_i-1 \}.
\end{split}
\end{equation*}

Set $A_a := r_{a!}\bar\Q_{l,[Y^a/G_{\nu\nu^a\nu'}]}
[\bar\nu_i+\nu_i-1]$. Identify $[\tE\nOi/\Gn]$ with the diagonal
part $\Delta$ of $[Y/\Gnn]$ and let $\iota: \Delta\to [Y/\Gnn]$ be
the inclusion. We have \setno
\begin{equation}\label{eqn:kef:iso62}
\begin{split}
  & \Id_{\D([\tE\nOi/\Gn])} \cong \pi'_! \iota_! \iota^* \pi^*
  \cong \pi'_!(\iota_!\bar\Q_{l,\Delta}\otimes \pi^*-), \\
  & \TE\Oi\TF\Oi = p_{1!} p'^*_1 p'_{1!} p_1^* [\bar\nu^1_i+\nu_i]
  \cong \pi'_! r_{1!} r_1^* \pi^* [\bar\nu_i+\nu_i-1]
  \cong \pi'_!(A_1\otimes \pi^*-), \\
  & \TF\Oi\TE\Oi = p_{2!} p'^*_2 p'_{2!} p_2^* [\nu^2_i+\bar\nu_i]
  \cong \pi'_! r_{2!} r_2^* \pi^* [\bar\nu_i+\nu_i-1]
  \cong \pi'_!(A_2\otimes \pi^*-).
\end{split}
\end{equation}

Note that $r_1,r_2$ are fibred respectively in $\BP^{\bar\nu_i-1}$
and $\BP^{\nu_i-1}$ over $\Delta$. Away from $\Delta$, both
$r_1,r_2$ restrict to an isomorphism onto the closed substack
$[U/\Gnn]\subset[Y/\Gnn]\setminus\Delta$ where $U\subset Y$ is such
that
\begin{align*}
  U/\Gnp\times GL(V_i) & = \{ (\dx,V,V') \mid \dim(V+V') = \nu_i+1 \} \\
  & = \{ (\dx,V,V') \mid \dim(V\cap V') = \nu_i-1 \}.
\end{align*}
It follows that
\begin{equation*}
\begin{split}
  & A_1|_\Delta \cong \bigoplus_{0\le m<\bar\nu_i} \bar\Q_{l,\Delta}[\bar\nu_i+\nu_i-1-2m], \\
  & A_2|_\Delta \cong \bigoplus_{0\le m<\nu_i} \bar\Q_{l,\Delta}[\bar\nu_i+\nu_i-1-2m], \\
  & A_1|_{[Y/\Gnn]\setminus\Delta} \cong A_2|_{[Y/\Gnn]\setminus\Delta}. \\
\end{split}
\end{equation*}

On the other hand, $[Y^a/G_{\nu\nu^a\nu'}]$ is a tower of
Grassmannian bundles over $[\tE\nOi/\Gn]$, hence is a smooth
algebraic stack and is proper over $[Y/\Gnn]$. By the decomposition
theorem, $A_1,A_2$ are semisimple complexes. Therefore, by the
classification of simple perverse sheaves, \setno
\begin{equation}\label{eqn:kef:iso63}
\begin{split}
  & A_1 \oplus \bigoplus_{\bar\nu_i\le m<\nu_i} \iota_!\bar\Q_{l,\Delta}[\bar\nu_i+\nu_i-1-2m] \\
  & \cong A_2 \oplus \bigoplus_{\nu_i\le m<\bar\nu_i} \iota_!\bar\Q_{l,\Delta}[\bar\nu_i+\nu_i-1-2m].
\end{split}
\end{equation}
From \eqref{eqn:kef:iso62}, \eqref{eqn:kef:iso63} the isomorphism
\eqref{eqn:kef:iso61} follows.

\subsubsection{Proof of (7)}

Let $\Omega\subset\hH$ be an orientation having $i$ as a source.
Assume $\nu'=\nu-i+j$. It suffices to show there is an isomorphism
of functors from $\D([\tE\nOi/\Gn])$ to $\D([\tE\npOi/\Gnp])$
\begin{equation*}
  \TE\Oi\TF_{\Omega,i,j} \cong \TF_{\Omega,i,j}\TE\Oi.
\end{equation*}

Let $\nu^1=\nu+j$, $\nu^2=\nu-i$. In a similar way as the proof of
(6), we reduce the problem to the consideration of the diagram,
$a=1,2$,
$$\xymatrixcolsep{1.5pc}\xymatrix{
  & [Y^a/G_{\nu\nu^a\nu'}] \ar[ld]_{q_a} \ar[rd]^{q'_a} \\
  [\tE\nOi/\Gn] & & [\tE\npOi/\Gnp] \\
}
$$
where $Y^a$ is the fibred product of some $\tZ\Oi$ over
$\tE_{\nu^a,\Omega,i}$ and
\begin{equation*}
\begin{split}
  & q_a(\dx,V,V',\dy) := (\dy^{-1}\dx\dy,V), \\
  & q'_a(\dx,V,V',\dy) := (\dx,V').
\end{split}
\end{equation*}
under the identification
\begin{equation*}
\begin{split}
  & Y^1/G_{\nu'}\times GL(V_i)\times GL(V^1_i)
  = \{ (\dx,V,V',\dy) \mid
  \dx\in E_{\nu^1,\dot\Omega}, \; \dy\in\dot F_{\nu\nu^1}, \;
  \\ & \quad\quad\quad\quad
  \Img \dx_h\dy_{h'}\subset\Img \dy_{h''}, \;
  V\subset V(i), \; V'\subset V^1=\dy(V), \;
  \\ & \quad\quad\quad\quad
  \dim V=\nu_i, \; \dim V'=\nu_i-1 \}, \\
  & Y^2/G_\nu\times GL(V^2_i)\times GL(V'_i)
  = \{ (\dx,V,V',\dy) \mid
  \dx\in E_{\nu',\dot\Omega}, \; \dy\in\dot F_{\nu^2\nu'}, \;
  \\ & \quad\quad\quad\quad
  \Img \dx_h\dy_{h'}\subset\Img \dy_{h''}, \;
  V^2\subset V\subset V(i), \; V'=\dy(V^2), \;
  \\ & \quad\quad\quad\quad
  \dim V=\nu_i, \; \dim V'=\nu_i-1 \}.
\end{split}
\end{equation*}
(Note that $E_{\nu,\dot\Omega}=E_{\nu^2,\dot\Omega}$,
$E_{\nu^1,\dot\Omega}=E_{\nu',\dot\Omega}$ and $\dot
F_{\nu\nu^1}=\dot F_{\nu^2\nu'}$.)

We have
\begin{equation*}
\begin{split}
  & \TE\Oi\TF_{\Omega,i,j}
  \cong q'_{1!} q_1^* [\bar\nu^1_j(\Omega)+\nu'_i]
  = q'_{1!} q_1^* [\nu_i+\bar\nu_j(\Omega)-2], \\
  & \TF_{\Omega,i,j}\TE\Oi
  \cong q'_{2!} q_2^* [\nu^2_i+\bar\nu'_j(\Omega)]
  = q'_{2!} q_2^* [\nu_i+\bar\nu_j(\Omega)-2].
\end{split}
\end{equation*}

Notice the isomorphism
$$s: [Y^2/G_{\nu\nu^2\nu'}] \to [Y^1/G_{\nu\nu^1\nu'}], \quad s(\dx,V,V',\dy) = (\dx,V,V',\dy).$$
We have $q_2=q_1s$ and $q'_2=q'_1s$. Hence $q'_{1!}q_1^* \cong
q'_{2!}q_2^*$. Hence $\TE\Oi\TF_{\Omega,i,j} \cong
\TF_{\Omega,i,j}\TE\Oi$.

\subsubsection{Proof of (9)}

Let $\Omega\subset\hH$ be an orientation having $i$ as a source.
Assume $\nu'=\nu+(1-a_{ij})i+j$. It suffices to show there is an
isomorphism of functors from $\D([\tE\nOi/\Gn])$ to
$\D([\tE\npOi/\Gnp])$
\begin{equation*}
  \bigoplus_{\stackrel{0\le m\le 1-a_{ij}}{\text{$m$ odd}}} \TF^{(1-a_{ij}-m)}\Oi\TF_{\Omega,i,j}\TF^{(m)}\Oi
  \cong \bigoplus_{\stackrel{0\le m\le 1-a_{ij}}{\text{$m$ even}}} \TF^{(1-a_{ij}-m)}\Oi\TF_{\Omega,i,j}\TF^{(m)}\Oi.
\end{equation*}

Let $\nu^m=\nu+mi$, $\nu'^m=\nu+mi+j$. In a similar way as the proof
of (6), we reduce the problem to the consideration of the diagram
$$\xymatrixcolsep{1.5pc}\xymatrix{
  & [Y^m/G_{\nu\nu^m\nu'^m\nu'}] \ar[d]^{r_m} \\
  [\tE\nOi/\Gn] & [Y/\Gnn] \ar[l]_--{\pi} \ar[r]^--{\pi'} & [\tE\npOi/\Gnp] \\
}
$$
where $Y,Y^m$ are algebraic varieties such that
\begin{equation*}
\begin{split}
  & Y/GL(V_i)\times GL(V'_i)
  = \{ (\dx,V,V',\dy) \mid \dx\in E_{\nu',\dot\Omega}, \; \dy\in\dot F_{\nu\nu'}, \;
  \\ & \quad\quad\quad\quad
  \Img \dx_h\dy_{h'}\subset\Img \dy_{h''}, \;
  V\subset V(i), \; \dy(V)\subset V'\subset V'(i), \;
  \\ & \quad\quad\quad\quad
  \dim V=\nu_i, \; \dim V'=\nu'_i \}, \\
  & Y^m/G_{\nu^m\nu'^m}\times GL(V_i)\times GL(V'_i)
  = \{ (\dx,V,V^m,V',\dy) \mid  \dx\in E_{\nu',\dot\Omega}, \; \dy\in\dot F_{\nu\nu'}, \;
  \\ & \quad\quad\quad\quad
  \Img \dx_h\dy_{h'}\subset\Img \dy_{h''}, \;
  V\subset V^m\subset V(i), \; \dy(V^m)\subset V'\subset V'(i), \;
  \\ & \quad\quad\quad\quad
  \dim V=\nu_i, \; \dim V^m=\nu^m_i, \; \dim V'=\nu'_i \}, \\
\end{split}
\end{equation*}
and
\begin{equation*}
\begin{split}
  & \pi(\dx,V,V',\dy) := (\dy^{-1}\dx\dy,V), \\
  & \pi'(\dx,V,V',\dy) := (\dx,V'), \\
  & r_m(\dx,V,V^m,V',\dy) := (\dx,V,V',\dy). \\
\end{split}
\end{equation*}

Set $A_m := r_{m!}\bar\Q_{l,[Y^m/G_{\nu\nu^m\nu'^m\nu'}]}$. We have
\begin{equation*}
  \TF^{(1-a_{ij}-m)}\Oi\TF_{\Omega,i,j}\TF^{(m)}\Oi
  \cong \pi'_! r_{m!} r_m^* \pi^* [t_m]
  \cong \pi'_! (A_m[t_m]\otimes\pi^*-)
\end{equation*}
where
\begin{equation*}
  t_m := m\bar\nu^m_i+\bar\nu'^m_j(\Omega)+(1-a_{ij}-m)\bar\nu'_i
  = \bar\nu'_j(\Omega)+(1-a_{ij})\bar\nu'_i+m(1-m).
\end{equation*}
So, we are left to show \setno
\begin{equation}\label{eqn:kef:iso9}
  \bigoplus_{\stackrel{0\le m\le 1-a_{ij}}{\text{$m$ odd}}} A_m[m(1-m)]
  \cong \bigoplus_{\stackrel{0\le m\le 1-a_{ij}}{\text{$m$ even}}} A_m[m(1-m)].
\end{equation}

Consider the stratification $Y=\bigsqcup_{n=1}^{1-a_{ij}}U_n$ where
$U_n$ is such that
$$U_n/GL(V_i)\times GL(V'_i) = \{ (\dx,V,V',\dy) \mid
  \dim(V'\cap\dy(V(i)))=\nu_i+n \}.
$$
Note that the morphism $r_m$ is fibred in the Grassmannian $Gr(m,n)$
of $m$-dimensional subspaces in an $n$-dimensional $\BF$-vector
space over $[U_n/\Gnn]$. So, it is a routine matter to verify that
the isomorphism \eqref{eqn:kef:iso9} holds on each stratum
$[U_n/\Gnn]$.

On the other hand, the variety $Y^m/G_{\nu^m\nu'^m}\times
GL(V_i)\times GL(V'_i)$ is a tower of Grassmannian bundles over
$$\{ (\dx,\dy) \in E_{\nu',\dot\Omega}\times \dot F_{\nu\nu'} \mid \Img \dx_h\dy_{h'}\subset\Img \dy_{h''} \}$$
and the latter is a vector bundle over $\dot F_{\nu\nu'}$. Hence
$[Y^m/G_{\nu\nu^m\nu'^m\nu'}]$ is a smooth algebraic stack and is
proper over $[Y/\Gnn]$. By the decomposition theorem, each $A_m$ is
a semisimple complex. Therefore, the validity of the isomorphism
\eqref{eqn:kef:iso9} follows from those on the strata $[U_n/\Gnn]$.

\subsection{Verdier duality and the bifunctor Ext}\label{sec:cat_fun:dual}

The propositions below show that Verdier duality and the bifunctor
$\Ext^\bullet_\FD(-,D-)$ categorify the bar involution of $U$ and
the contravariant form of $U$-module, respectively.

\begin{prop}\label{prop:kef:dual}
We have isomorphisms of endofunctors of $\FD$
\begin{equation*}
  D[-1]=[1]D, \quad
  D\FK_i=\FK_i^{-1}D, \quad
  D\FE^{(n)}_i=\FE^{(n)}_iD, \quad
  D\FF^{(n)}_i=\FF^{(n)}_iD.
\end{equation*}
\end{prop}

The first two isomorphisms are obvious.

Recall from \ref{sec:cat_fun:alt:p} that the morphisms $p_i,p'_i$
are proper and smooth of relative dimension $n\bar\nu'_i$, $n\nu_i$,
respectively. Hence $Dp_{i!}=p_{i!}D$ and
$Dp'^*_i[n\nu_i]=p'^*_i[n\nu_i]D$. It follows that
$D\TE^{(n)}\Oi=\TE^{(n)}\Oi D$. Similarly,
$D\TF^{(n)}\Oi=\TF^{(n)}\Oi D$. These conclude the last two
isomorphisms (cf. \ref{sec:pre:loc:transform}).

\begin{prop}\label{prop:kef:form}
There are natural isomorphisms for $A,B\in\FD$
\begin{equation*}
\begin{array}{l@{}ll}
  \Ext^\bullet_\FD(A,DB) & = \Ext^\bullet_\FD(B,DA), \\
  \Ext^\bullet_\FD(\FK_iA,DB) & = \Ext^\bullet_\FD(A,D\FK_iB), \\
  \Ext^\bullet_\FD(\FE_i^{(n)}A,DB) & = \Ext^\bullet_\FD(A,D\FK_i^n\FF_i^{(n)}[-n^2]B), \\
  \Ext^\bullet_\FD(\FF_i^{(n)}A,DB) & = \Ext^\bullet_\FD(A,D\FK_i^{-n}\FE_i^{(n)}[-n^2]B).
\end{array}
\end{equation*}
\end{prop}

The first one actually states that Verdier duality $D$ is self
adjoint, which is an obvious fact.

The others follow from \propref{prop:kef:adjoint} and
\propref{prop:kef:dual}.

\section{Categorification of integrable representations}\label{sec:rep}

Let the quiver $(I,H)$, the generalized Cartan matrix $(a_{ij})$ and
the quantum group $U$ be defined in \ref{sec:pre:u}.

\setcounter{subsubsection}{0}

\subsubsection{Notations}

An element $\mu=\sum_{i\in\hI}\mu_i\cdot i\in\N[\hI]$ is called of
\begin{enumerate}
  \item type I, if $\mu=ni$ for some $n\ge1$ and $i\in I$;
  \item type II, if $\mu\in\N[\hI\setminus I]$, i.e. $\mu_i=0$ for $i\in I$.
\end{enumerate}

In this section, we use $\omega$ (resp. $\vo,\vm$) to denote a type
II element (resp. a sequence of type II elements, a sequence of type
I,II elements) in $\N[\hI]$. We write $\vm\rhd\vo$ for the fact that
the subsequence of $\vm$ formed by the type II terms is identical to
$\vo$.

Let $M(\omega),R(\omega),\Lambda(\omega),\Lambda(\vo)$ be the
$U$-modules defined in \ref{sec:pre:u:module} by regarding
$\N[\hI\setminus I]$ as $\N[I]$.

For a sequence $\vn=(\nu^1,\nu^2,\dots,\nu^s)$ we set
$|\vn|=\sum_{a=1}^s\nu^a$ and write $\vn\mu$ for the sequence
$(\nu^1,\nu^2,\dots,\nu^s,\mu)$.

\subsection{The functor $\FF_\mu$}\label{sec:rep:f}

First, we fix a full flag
$$\BF^0\subset\BF^1\subset\BF^2\subset\BF^3\subset\cdots$$
so that there is an unambiguous inclusion $\BF^m\subset\BF^n$ for
$m\le n$.

\subsubsection{}

Suppose $\mu=\sum_{i\in\hI}\mu_i\cdot i\in\N[\hI]$ has discreet
support, i.e. $\mu_{h'}\mu_{h''}$ vanishes for all $h\in\hH$. In
particular, $\mu$ can be any type I or type II element.

Assume $\nu'=\nu+\mu$. Using the variety
$$F_{\nu\nu'} := \{ y\in\bigoplus_{j\in\hI}\Hom(V_j,V'_j) \mid \Ker y_j=0, \;
  \text{$y_{\hj}$ is the inclusion,} \; j\in I \}
$$
we generalize the functor $\F^{(n)}\Oi$ defined in
\ref{sec:cat_fun:fun:kef} to the following
$$\F_{\Omega,\mu} := p'_!p^*[\bar\nu'_\mu(\Omega)] : \D([E\nO/\Gn])\to\D([E\npO/\Gnp]).$$
where
$$\bar\nu_\mu(\Omega) := \sum_{i\in I} \mu_i\bar\nu_i(\Omega)
  + \sum_{i\in\hI\setminus I} \mu_i \sum_{h\in\Omega: \; h'=i} \nu_{h''}.
$$
Note that $\F^{(n)}\Oi = \F_{\Omega,ni}$.

\begin{prop}
We have $\Phi_{\Omega_1,\Omega_2}\F_{\Omega_1,\mu} \cong
\F_{\Omega_2,\mu}\Phi_{\Omega_1,\Omega_2}$ and
$\F_{\Omega,\mu}(\CN\nO)\subset\CN\npO$ for orientations
$\Omega,\Omega_1,\Omega_2\subset\hH$.
\end{prop}

The proof is the same as \propref{prop:kef:fourier} and the easy
part of \propref{prop:kef:local}.

\subsubsection{}

It follows that $\F_{\Omega,\mu}$ induces a well defined functor
$$\FF_{\nu,\mu} : \FD_\nu\to\FD_{\nu+\mu}.$$
Define an endofunctor of $\FD$ as usual
$$\FF_\mu := \bigoplus_\nu \FF_{\nu,\mu}.$$

\subsubsection{}\label{sec:rep:f:iota}

Suppose $\mu$ has discreet support. Let $\Omega\subset\hH$ be an
orientation having the support of $\mu$ as sinks, i.e. $i\in\hI$ is
a sink of $\Omega$ for all nonvanishing $\mu_i$.

Assume $\nu=\nu'+\mu$. We have an unambiguous inclusion
$E\npO\subset E\nO$. Moreover, regarding $\Gnp$ as a subgroup of
$\Gn$ so that the inclusion $E\npO\subset E\nO$ is
$\Gnp$-equivariant, we get a representable morphisms $\iota_\mu:
[E\npO/\Gnp] \to [E\nO/\Gn]$ which is unique up to isomorphism.

\begin{prop}\label{prop:rep:f}
Let $\Omega$ be an orientation and let $\mu$ be of type II.
\begin{enumerate}
\item
$\FK_i\FF_\mu = \FF_\mu\FK_i[-\mu_{\hi}]$ and $\FE^{(n)}_i\FF_\mu
\cong \FF_\mu\FE^{(n)}_i$;
\item
$D\F_{\Omega,\mu} \cong \F_{\Omega,\mu} D$, $D\FF_\mu \cong \FF_\mu
D$;
\item
$\F_{\Omega,\mu}$ and $\FF_\mu$ are fully faithful;
\item
$\F_{\Omega,\mu}$ and $\FF_\mu$ send simple perverse sheaves to
simple perverse sheaves.
\end{enumerate}
\end{prop}

The rest of this subsection is dedicated to the proof of the
proposition.

Assume $\Omega$ has the vertices $\hI\setminus I$ as sinks and
$\nu=\nu'+\mu$. Note that $\Gn=\Gnp$ thus $\iota_\mu: [E\npO/\Gnp]
\to [E\nO/\Gn]$ is an inclusion of closed substack.

\begin{lem}
$\F_{\Omega,\mu}=\iota_{\mu!}$.
\end{lem}

Our assumption on $\Omega$ implies $\bar\nu'_\mu(\Omega)=0$. Note
that $p: [Z_\Omega/\Gnn]\to[E\nO]$ is an isomorphism and
$p'=\iota_\mu p$. Therefore,
$\F_{\Omega,\mu}=p'_!p^*=\iota_{\mu!}p_!p^*=\iota_{\mu!}$.

\begin{lem}\label{lem:rep:f2}
Let $\Omega,\Omega_1,\Omega_2\subset\hH$ be orientations all having
the vertices $\hI\setminus I$ as sinks. We have
$\Phi_{\Omega_1,\Omega_2}\iota_\mu^* \cong
\iota_\mu^*\Phi_{\Omega_1,\Omega_2}$ and
$\iota_\mu^*(\CN\nO)\subset\CN\npO$.
\end{lem}

The proof is the same as \propref{prop:kef:fourier} and the easy
part of \propref{prop:kef:local}.

\subsubsection{Proof of \propref{prop:rep:f}}

Assume $\Omega$ has the vertices $\hI\setminus I$ as sinks.

(1) The first isomorphism is obvious. To see the second one, we
assume $i$ is a source of $\Omega$ and form the following cartesian
squares for $\tilde\nu=\nu+ni$, $\tilde\nu'=\nu'+ni$
\begin{equation*}
\xymatrixcolsep{1.5pc}\xymatrix{
  [E\nO/\Gn] \ar[d]_{\iota_\mu}
  & [\tZ\Oi/G_{\nu\tilde\nu}] \ar[l]_--{\tp_i} \ar[r]^--{\tp'_i} \ar[d]
  & [E_{\tilde\nu,\Omega}/G_{\tilde\nu}] \ar[d]^{\iota_\mu} \\
  [E\npO/\Gnp]
  & [\tZ'\Oi/G_{\nu'\tilde\nu'}] \ar[l]_--{\tp_i} \ar[r]^--{\tp'_i}
  & [E_{\tilde\nu',\Omega}/G_{\tilde\nu'}]
}
\end{equation*}
By proper base change
$$\E^{(n)}_{\Omega,i}\F_{\Omega,\mu}
  = \tp_{i!}\tp'^*_i[n\nu'_i]\iota_{\mu!}
  \cong \iota_{\mu!}\tp_{i!}\tp'^*_i[n\nu_i]
  = \F_{\Omega,\mu}\E^{(n)}_{\Omega,i}.
$$
The second isomorphism follows.

(2) is immediate from the isomorphism
$D\iota_{\mu!}\cong\iota_{\mu!}D$.

(3) Clearly $\F_{\Omega,\mu}=\iota_{\mu!}$ is fully faithful. By
\lemref{lem:rep:f2}, $\iota_\mu^*$ gives rise to an endofunctor
$\FRes_\mu$ of $\FD$. Since $\iota_\mu^*$ is left adjoint to
$\iota_{\mu!}=\iota_{\mu*}$ and since the adjunction morphism
$\iota_\mu^*\iota_{\mu*}\to\Id$ is an isomorphism, $\FRes_\mu$ is
left adjoint to $\FF_\mu$ and the adjunction morphism
$\FRes_\mu\FF_\mu\to\Id$ is an isomorphism. Therefore, the functor
$\FF_\mu$ is fully faithful.

(4) Clearly $\F_{\Omega,\mu}=\iota_{\mu!}$ sends simple perverse
sheaves to simple perverse sheaves. Thus $\FF_\mu$ sends a simple
perverse sheave either to a simple perverse sheaf or to zero (cf.
\ref{sec:pre:loc:t}). But the latter case may not happen, for
$\FF_\mu$ is fully faithful.

\subsection{The subcategory $\FQ_\vo$ of $\FD$}\label{sec:rep:q}

\subsubsection{}\label{sec:rep:q:l}

Note that the variety $E_{0,\Omega}$ is a single point and $G_0$ is
a trivial group. Thus $\FD_0 = \D([E_{0,\Omega}/G_0])/\CN_{0,\Omega}
= \D([E_{0,\Omega}/G_0]) = \D(\pt)$ for every orientation
$\Omega\subset\hH$.

So, for a sequence $\vm=(\mu^1,\mu^2,\dots,\mu^s)$ of type I,II
elements in $\N[\hI]$, we have a well defined complex
$$\CL_\vm := \FF_{\mu^s} \cdots \FF_{\mu^2} \FF_{\mu^1} \bar\Q_{l,\pt} \in \FD.$$
By definition $\CL_\vm$ is represented by
$$\CL_{\Omega,\vm} = \F_{\Omega,\mu^s} \cdots \F_{\Omega,\mu^2} \F_{\Omega,\mu^1} \bar\Q_{l,\pt}.$$

\begin{prop}\label{prop:rep:L}
$\CL_{\Omega,\vm}$ and hence $\CL_\vm$ are semisimple complexes.
\end{prop}

Assume $\vm=(\mu^1,\mu^2,\dots,\mu^s)$. Set $\nu := |\vm|$, $|\vm|^a
:= \sum_{b\le a} \mu^b$, $F_\vm :=
\prod_{a=1}^sF_{|\vm|^{a-1}|\vm|^a}$, $G_\vm :=
\prod_{a=1}^sG_{|\vm|^a}$ and
\begin{equation*}
  Z_\vm := \{ (x,y) \in E\nO\times F_\vm \mid
  \Img x_h(y_sy_{s-1}\cdots y_a)_{h'}\subset\Img (y_sy_{s-1}\cdots y_a)_{h''} \}.
\end{equation*}
Let $\pi_\vm: [Z_\vm/G_\vm] \to [E\nO/\Gn]$ be the obvious
projection. An easy proper base change argument shows that
$$\CL_{\Omega,\vm} \cong \pi_{\vm!}\bar\Q_{l,[Z_\vm/G_\vm]} [\vm(\Omega)]$$
where
$$\vm(\Omega) := \sum_{a<b}(\sum_{h\in\Omega}\mu^a_{h''}\mu^b_{h'}-\sum_{i\in I}\mu^a_i\mu^b_i).$$

Let $G:=\prod_{a=1}^{s-1}G_{|\vm|^a}$. Note that $Z_\vm/G$ is a
vector bundle over the partial flag variety $F_\vm/G$ and is proper
over $E\nO$. Thus $[Z_\vm/G_\vm]$ is smooth and $\pi_\vm$ is
representable and proper. By the decomposition theorem
$\CL_{\Omega,\vm}$ is a semisimple complex.

\subsubsection{}

Let $\FQ_\vo$ be the full subcategory of $\FD$ formed by the finite
direct sums $\oplus_r A_r[n_r]$ where $n_r\in\Z$ and $A_r$ is a
direct summand of some $\CL_\vm$ with $\vm\rhd\vo$.

Similarly define a full subcategory $\CQ_{\Omega,\vo}$ of
$\bigoplus_\nu\D([E\nO/\Gn])$ for an orientation $\Omega$.

\begin{prop}\label{prop:rep:kef}
$\FQ_\vo$ is stable under the functors $\FK^{\pm1}_i, \FE^{(n)}_i,
\FF^{(n)}_i$.
\end{prop}

The claim is clear for $\FK^{\pm1}_i, \FF^{(n)}_i$.

By \thmref{thm:kef}(4) $\FE^{(n)}_iA$ is a direct summand of
$\FE^n_iA[\frac{n(n-1)}2]$ for $A\in\FD$. So it suffices to show
$\FE_i\CL_\vm\in\FQ_\vo$ for $\vm\rhd\vo$. By a similar reason, we
may assume the type I terms in $\vm$ are in the form $j$ for various
$j\in I$.

We prove by induction on the length of $\vm$.

If $\vm=\emptyset$, then $\FE_i\CL_\vm=\FE_i\bar\Q_{l,\pt}=0$; we
are done.

If $\vm=\vm'\omega$, by \propref{prop:rep:f}(1)
$\FE_i\CL_\vm=\FE_i\FF_\omega\CL_{\vm'}=\FF_\omega\FE_i\CL_{\vm'}$,
which by the inductive hypothesis is contained in $\FQ_\vo$.

If $\vm=\vm'j$ for some $j\in I$, by \thmref{thm:kef}(6)(7)
$\FE_i\CL_\vm=\FE_i\FF_j\CL_{\vm'}$ is either (i) isomorphic to
$\FF_j\FE_i\CL_{\vm'}$ or (ii) isomorphic to a direct summand of
$\FF_j\FE_i\CL_{\vm'}$ or (iii) isomorphic to the direct sum of
$\FF_j\FE_i\CL_{\vm'}$ and several $\CL_{\vm'}[n_r]$. In any case,
$\FE_i\CL_\vm$ is contained in $\FQ_\vo$ by the inductive
hypothesis.

\begin{prop}\label{prop:rep:form}
$\sum_n \dim\Ext^n_\FD(A,B)\cdot q^{-n} \in \N[q,q^{-1}]$ for
$A,B\in\FQ_\vo$.
\end{prop}

Assume $A,B\in\FQ_\vo\cap\FD_\nu$ and assume both are in the form
$\CL_\vm$ with $\vm\rhd\vo$. We prove by induction on the norm
$\|\nu\|:=\sum_{i\in\hI}\nu_i$ and the length of $\vo$.

If $\nu=0$, then $A=B=\bar\Q_{l,\pt}$; we are done.

If either of $A,B$, say $A$, has the form $\CL_{\vm ni}$ with
$\vm\rhd\vo$, $n\ge1$, $i\in I$, then by \propref{prop:kef:adjoint}
$$\Ext^\bullet_\FD(A,B)
  = \Ext^\bullet_\FD(\FF_i^{(n)}\CL_\vm,B)
  = \Ext^\bullet_\FD(\CL_\vm,\FK_i^n\FE_i^{(n)}[n^2]B).
$$
Thus by the inductive hypothesis, the proposition is true for $A,B$.

If $A=\CL_{\vm\omega}, B=\CL_{\vm'\omega}$ such that
$\vm\omega,\vm'\omega\rhd\vo$, then by \propref{prop:rep:f}(3)
$$\Ext^\bullet_\FD(A,B)
  = \Ext^\bullet_\FD(\FF_\omega\CL_\vm,\FF_\omega\CL_{\vm'})
  = \Ext^\bullet_\FD(\CL_\vm,\CL_{\vm'}).
$$
Thus by the inductive hypothesis, the proposition is also true for
$A,B$.

\begin{prop}\label{prop:rep:sink}
Let $A\in\CQ_{\Omega,\vo}$ be a simple perverse sheaf. Then
\begin{enumerate}
\item
$A$ is self dual, i.e. $DA\cong A$ (therefore, both
$\CQ_{\Omega,\vo}$ and $\FQ_\vo$ are stable under Verdier duality);
and
\item
there is an isomorphism with $\vm_r,\vm'_r\rhd\vo$
$$A \oplus \CL_{\Omega,\vm_1}[n_1] \oplus \CL_{\Omega,\vm_2}[n_2] \oplus \cdots \oplus \CL_{\Omega,\vm_s}[n_s]
  \cong \CL_{\Omega,\vm'_1}[n'_1] \oplus \CL_{\Omega,\vm'_2}[n'_2] \oplus \cdots \oplus \CL_{\Omega,\vm'_t}[n'_t].
$$
\end{enumerate}
\end{prop}

We prove the proposition in the rest of this subsection.

\subsubsection{}

Let $\Omega\subset\hH$ be an orientation having a vertex $i\in I$ as
a sink. We denote by $\bar x(i)$ the restriction of $x\in E\nO$ to
the direct summand
$$\bigoplus_{h\in\Omega:\;h''=i}\Hom(V_{h'},V_{h''}) = \Hom(\bigoplus_{h\in\Omega:\;h''=i}V_{h'},V_i).$$
There is a filtration of closed substacks
$$[E\nO/\Gn]=U_{\ge0} \supset U_{\ge1} \supset \cdots \supset U_{\ge\nu_i} \supset U_{\ge\nu_i+1}=\emptyset$$
where
$$U_{\ge n} := [\{ x\in E\nO \mid \dim\Coker \bar{x}(i) \ge n \}/\Gn].$$
Let $U_n$ be the stratum
$$U_n := U_{\ge n}\setminus U_{\ge n+1} = [\{ x\in E\nO \mid \dim\Coker \bar{x}(i) = n \}/\Gn].$$

\begin{lem}\label{lem:rep:sink1}
Assume $\nu'=\nu+ni$ for some $i\in I$. Then
$Supp(\F_{\Omega,ni}A)\subset U'_{\ge n}$ for $A\in\D([E\nO/\Gn])$.
\end{lem}

This is clear from the inclusion $q'([Z_\Omega/\Gnn]) \subset
U'_{\ge n}$.

\subsubsection{}

Assume $\nu=\nu'+ni$. Let $\iota_n: [E\npO/\Gnp] \to [E\nO/\Gn]$ be
the representable morphism induced by the inclusion $E\npO\subset
E\nO$ (cf. \ref{sec:rep:f:iota}).

Note that $\iota_n$ restricts to a representable morphism
$\iota_n|_{U'_0}: U'_0\to U_n$ which is smooth with connected
nonempty fibers of dimension $n\nu_i$. Thus the functor
$(\iota_n|_{U'_0})^*[n\nu_i]$ is perverse t-exact, commutes with
Verdier duality and induces a fully faithful functor
$\M(U_n)\to\M(U'_0)$.

\begin{lem}\label{lem:rep:sink2}
Let $A\in\D([E\nO/\Gn])$ be a simple perverse sheaf. Assume
$Supp(A)\subset U_{\ge n}$ and $Supp(A)\cap U_n\ne\emptyset$. Then
\begin{enumerate}
\item
$\iota_n^*[n\nu_i]A|_{U'_0}$ is a simple perverse sheaf.
\end{enumerate}
Moreover, let $B\in\D([E\npO/\Gnp])$ be the intermediate extension
of $\iota_n^*[n\nu_i]A|_{U'_0}$.
\begin{enumerate}
\setcounter{enumi}{1}
\item
If $B$ is self dual, so is $A$.
\item
$A|_{U_n}\cong\F_{\Omega,ni}B|_{U_n}$.
\end{enumerate}
\end{lem}

First, our assumptions on $A$ imply that $A$ is the intermediate
extension of the simple perverse sheaf $A|_{U_n}$. It follows that
$\iota_n^*[n\nu_i]A|_{U'_0}=(\iota_n|_{U'_0})^*[n\nu_i](A|_{U_n})$
is a simple perverse sheaf. This proves (1).

Moreover, if $B$ is self dual, so is
$B|_{U'_0}=\iota_n^*[n\nu_i]A|_{U'_0}=(\iota_n|_{U'_0})^*[n\nu_i](A|_{U_n})$,
hence so is $A|_{U_n}$, hence so is $A$. This proves (2).

Next, we form the following commutative diagram
\begin{equation*}
\xymatrixcolsep{1.5pc}\xymatrix{
  [Z_\Omega/\Gnn] \ar[rd]_p & E\npO \ar[l]_--{f} \ar[r]^--{f} \ar[d]^{\pi} & [Z_\Omega/\Gnn] \ar[d]^{p'} \\
  & [E\npO/\Gnp] \ar[r]_{\iota_n} & [E\nO/\Gn]
}
\end{equation*}
where $Z_\Omega,p,p'$ are those defining $\F_{\Omega,ni}$, $\pi$ is
the presentation of the quotient stack and $f(x):=(\iota
x\iota^{-1},\iota)$, $\iota\in F_{\nu'\nu}$ being the inclusion.
($f$ is well defined since $i$ is a sink of $\Omega$.)

The commutative diagram restricts to the following one
\begin{equation*}
\xymatrixcolsep{3pc}\xymatrix{
  U_Z \ar[rd]_{p|_{U_Z}} & \pi^{-1}(U'_0) \ar[l]_{f} \ar[r]^-{f} \ar[d]^{\pi} & U_Z \ar[d]^{p'|_{U_Z}} \\
  & U'_0 \ar[r]_{\iota_n|_{U'_0}} & U_n
}
\end{equation*}
where
\begin{equation*}
  U_Z := [\{ (x,y)\in Z_\Omega \mid \dim\Coker \bar{x}(i) = n \}/\Gnn].
\end{equation*}

Note that $p'|_{U_Z}: U_Z\to U_n$ is an isomorphism. We have
\begin{equation*}
\begin{split}
  & \pi^*(\iota_n|_{U'_0})^*[n\nu_i](\F_{\Omega,ni}B|_{U_n})
  = \pi^*(\iota_n|_{U'_0})^*(p'_!p^*B|_{U_n}) \\
  & = f^*(p'|_{U_Z})^*(p'|_{U_Z})_!(p|_{U_Z})^*(B|_{U'_0}) \\
  & \cong \pi^*(B|_{U'_0})
  = \pi^*(\iota_n|_{U'_0})^*[n\nu_i](A|_{U_n}).
\end{split}
\end{equation*}
Since the functors $\pi^*[\dim\Gnp]$ and
$(\iota_n|_{U'_0})^*[n\nu_i]$ are perverse t-exact and induce fully
faithful functors $\M(U_n)\to\M(U'_0)\to\M(\pi^{-1}(U'_0))$, it
follows that $A|_{U_n}\cong\F_{\Omega,ni}B|_{U_n}$. This proves (3).

\begin{lem}\label{lem:rep:sink3}
Assume $\vm=(\mu^1,\mu^2,\dots,\mu^s)$ and $\nu=\nu'+ni=|\vm|$. We
have
$$\iota_n^*\CL_{\Omega,\vm} \cong \bigoplus_{\vm'\in M} \CL_{\Omega,\vm'}
  \boxtimes u_!\bar\Q_{l,P_{\vm'}}[n_{\vm'}]
$$
for some integers $n_{\vm'}$, where
$$M := \{ \vm' \mid \mu'^a_j\le\mu^a_j, \; |\vm'|=|\vm|-ni \},$$
and $u: P_{\vm'}\to\pt$ is the point map of the partial flag variety
$$P_{\vm'} := \{ (0=V_0\subset V_1\subset V_2\subset\cdots\subset V_s) \mid
  \dim V_a/V_{a-1}=\mu^a_i-\mu'^a_i \}.
$$
\end{lem}

Keep the notations of \propref{prop:rep:L}. We set
\begin{equation*}
\begin{split}
  & Z := \{ (x,y)\in Z_\vm \mid x\in E\npO \}, \\
  & G:=\prod_{a=1}^{s-1}G_{|\vm|^a}\times\Gnp.
\end{split}
\end{equation*}
Let $\pi: [Z/G]\to [E\npO/\Gnp]$ be the obvious projection. By
proper base change
$$\iota_n^*\CL_{\Omega,\vm} \cong \pi_!\bar\Q_{l,[Z/G]}.$$

We have a stratification $[Z/G]=\sqcup_{\vm'\in M}[U_{\vm'}/G]$
where
$$U_{\vm'} :=  \{ (x,y)\in Z \mid \dim (\Img(y_sy_{s-1}\cdots y_a)_i\cap V'_i) = |\vm'|^a_i \}.$$
Note that restricting to the stratum $[U_{\vm'}/G]$, $\pi$ can be
factored as
$$[U_{\vm'}/G] \xrightarrow{f_{\vm'}} Z_{\vm'}\times P_{\vm'} \xrightarrow{\pi_{\vm'}\times u} [E\npO/\Gnp]$$
where $f_{\vm'}$ is a vector bundle of fiber dimension
$\sum_{a<b}(\mu^a_i-\mu'^a_i)\mu'^b_i$. Therefore, from the
decomposition theorem \ref{sec:pre:perverse:decom} our lemma
follows.

\subsubsection{Proof of \propref{prop:rep:sink}}

Assume
$A\in\CQ_{\Omega,\vo,\nu}:=\CQ_{\Omega,\vo}\cap\D([E\nO/\Gn])$. We
show the proposition by induction on the norm
$\|\nu\|:=\sum_{i\in\hI}\nu_i$ and the length of $\vo$.

Case (1) $\nu=0$. That is,
$A=\CL_{\Omega,\emptyset}=\bar\Q_{l,\pt}$. We are done.

Case (2) $A$ is a direct summand of
$\F_{\Omega,\omega}\CL_{\Omega,\vm}[m]$ with $\vm\omega\rhd\vo$. By
\propref{prop:rep:f}(4), $A\cong\F_{\Omega,\omega}B$ where $B$ is a
direct summand of $\CL_{\Omega,\vm}[m]$. Then by
\propref{prop:rep:f}(2) and by the inductive hypothesis, both claims
of the proposition are true for $A$.

Case (3) $A$ is a direct summand of
$\F_{\Omega,mi}\CL_{\Omega,\vm}[m']$ with $\vm\rhd\vo$, $m\ge1$,
$i\in I$. Assume $i$ is a sink of $\Omega$ and assume
$Supp(A)\subset U_{\ge n}$ and $Supp(A)\cap U_n\ne\emptyset$. By
\lemref{lem:rep:sink1}, $n\ge m$.

Let $\nu'=\nu-ni$. By \lemref{lem:rep:sink3}, $\iota_n^*A$ is a
semisimple complex in $\CQ_{\Omega,\vo,\nu'}$. By
\lemref{lem:rep:sink2}(1), $\iota_n^*[n\nu_i]A|_{U'_0}$ is a simple
perverse sheaf. Hence the intermediate extension $B$ of
$\iota_n^*[n\nu_i]A|_{U'_0}$ is a simple perverse sheaf in
$\CQ_{\Omega,\vo,\nu'}$, which is self dual by the inductive
hypothesis. Then by \lemref{lem:rep:sink2}(2)(3), $A$ is self dual
and $\F_{\Omega,ni}B\cong A\oplus C$ with $C\in\CQ_{\Omega,\vo,\nu}$
and $Supp(C)\subset U_{\ge n+1}$.

Applying the argument in the last paragraph on the simple direct
summands of $C$, and so on, yields
$$A \oplus A_1[m_1] \oplus A_2[m_2] \oplus \cdots \oplus A_s[m_s]
  \cong A'_1[m'_1] \oplus A'_2[m'_2] \oplus \cdots \oplus A'_t[m'_t]
$$
where $A_r,A'_r$ are in the form $\F_{\Omega,ni}B$ with $1\le n\le
\nu_i$ and $B$ being a simple perverse sheaf in
$\CQ_{\Omega,\vo,\nu-ni}$. By the inductive hypothesis, the second
claim of the proposition holds for $A_r,A_r'$, thus holds for $A$.

\subsection{Main theorems}\label{sec:rep:main}

\subsubsection{}

Let $\bG(\FQ_\vo)$ denote the Grothendieck group of the additive
category $\FQ_\vo$. That is, $\bG(\FQ_\vo)$ is the free
$\Z[q,q^{-1}]$-module defined by the generators each for an
isomorphism class of objects from $\FQ_\vo$ and the relations
\begin{enumerate}
\renewcommand{\theenumi}{\roman{enumi}}
  \item $[A\oplus B]=[A]+[B]$, for $A,B\in\FQ_\vo$;
  \item $[A[1]]=q^{-1}[A]$, for $A\in\FQ_\vo$.
\end{enumerate}

Set $\tG(\FQ_\vo):=\bG(\FQ_\vo)\otimes\Q(q)$. By definition it has a
basis
$$\CB(\vo) := \{ \; [A] \mid \text{$A\in\FQ_\vo$ is a simple perverse sheaf in $\FD$} \}.$$
Moreover, the functor $\FF_\omega$ induces a linear homomorphism
$$\varphi_\omega: \tG(\FQ_\vo) \to \tG(\FQ_{\vo\omega}).$$

\begin{thm}
The followings endow $\tG(\FQ_\vo)$ with a structure of
$U$-module.
$$K^{\pm1}_i[A]:=[\FK^{\pm1}_iA], \quad E^{(n)}_i[A]:=[\FE^{(n)}_iA], \quad F^{(n)}_i[A]:=[\FF^{(n)}_iA].$$
\end{thm}

This is the consequence of \thmref{thm:kef} and
\propref{prop:rep:kef}.

\begin{thm}\label{thm:rep:form}
The following defines a nondegenerate contravariant form on the
$U$-module $\tG(\FQ_\vo)$.
$$([A],[B]) := \sum_n \dim\Ext^n_\FD(A,DB)\cdot q^{-n}.$$
\end{thm}

By \propref{prop:rep:form} and \propref{prop:rep:sink}(1), the above
expression values in $\N[q,q^{-1}]$ for $A,B\in\FQ_\vo$. So, by
\propref{prop:kef:form} it defines a contravariant form of
$U$-module.

It is a basic property of a t-structure that for simple perverse
sheaves $A,B\in\FD$, $\Ext^n_\FD(A,B)$ vanishes for $n<0$, and
$\Ext^0_\FD(A,B)\cong\bar\Q_l$ if $A\cong B$ or vanishes otherwise.
By \propref{prop:rep:sink}(1), simple perverse sheaves in $\FQ_\vo$
are self dual. Hence $(b,b')\in\delta_{bb'}+q^{-1}\N[q^{-1}]$ for
$b,b'\in\CB(\vo)$. Therefore, the contravariant form under the basis
$\CB(\vo)$ yields a unit matrix modulo $q^{-1}$. This proves the
nondegeneracy.

\begin{thm}
There exist a unique family of isomorphisms of $U$-modules
$\Lambda(\vo) \cong \tG(\FQ_\vo)$ so that
\begin{enumerate}
  \item the isomorphisms preserve contravariant form,
  \item identifying $\Lambda(\vo)$ with $\tG(\FQ_\vo)$, we have $\CB(\emptyset) = \{1\}$ and
  \item $\varphi_\omega(u)=u\otimes\eta_\omega$ for $u\in\tG(\FQ_\vo)$.
\end{enumerate}
\end{thm}

By \propref{prop:rep:f}(1) $E_i\varphi_\omega(u) =
\varphi_\omega(E_iu)$ and $K_i\varphi_\omega(u) =
q^{\omega_{\hi}}\varphi_\omega(K_iu)$ for $u\in\tG(\FQ_\vo)$.
Moreover, by \propref{prop:rep:f}(2)(3), the homomorphism
$\varphi_\omega$ preserves contravariant form.

Therefore, by the universality of $M(\omega)$ (cf.
\ref{sec:pre:u:verma}) there exists a unique homomorphism of
$U$-modules
$$\tilde\varphi_\omega: \tG(\FQ_\vo)\otimes M(\omega) \to \tG(\FQ_{\vo\omega})$$
such that
$\tilde\varphi_\omega(u\otimes\eta_\omega)=\varphi_\omega(u)$,
$u\in\tG(\FQ_\vo)$. Moreover, $\tilde\varphi_\omega$ preserves
contravariant form.

Since the contravariant forms on $\tG(\FQ_\vo),\tG(\FQ_{\vo\omega})$
are nondegenerate, the kernel of $\tilde\varphi_\omega$ is
$\tG(\FQ_\vo)\otimes R(\omega)$. On the other hand,
\propref{prop:rep:sink}(2) implies that the set $\{ \, [\CL_\vm]
\mid \vm\rhd\vo\omega \}$ generates $\tG(\FQ_{\vo\omega})$ as a
linear space.  Thus $\varphi_\omega(\tG(\FQ_\vo))$ generates
$\tG(\FQ_{\vo\omega})$ as a $U$-module. Thus $\tilde\varphi_\omega$
is an epimorphism. Therefore $\tilde\varphi_\omega$ induces an
isomorphism of $U$-modules $\tG(\FQ_\vo)\otimes \Lambda(\omega)
\cong \tG(\FQ_{\vo\omega})$.

Now observe that $\F_{\Omega,ni}\bar\Q_{l,\pt} \in
\CN_{ni,\Omega,i}$ for $n\ge1$, $i\in I$. Thus $\CL_\vm$ vanishes
for nonempty $\vm\rhd\emptyset$. Thus $\FQ_\emptyset\sim\D(\pt)$ and
$\tG(\FQ_\emptyset)$ is a trivial $U$-module. So, there is a unique
isomorphism of $U$-modules
$\Lambda(\emptyset)\cong\tG(\FQ_\emptyset)$ sending $1$ to the
unique element $[\bar\Q_{l,\pt}]$ of $\CB(\emptyset)$.

Then an induction on the length of $\vo$ establishes the desired
isomorphisms $\Lambda(\vo) \cong \tG(\FQ_{\vo})$. The uniqueness
is clear from the proof.

\begin{thm}\label{thm:rep:bar}
Verdier duality induces a family of $\Q$-linear involutions $\Psi:
\tG(\FQ_\vo) \to \tG(\FQ_\vo)$. They satisfy
\begin{enumerate}
  \item $\Psi(xu) = \bar{x}\Psi(u)$, for $x\in U$, $u \in \tG(\FQ_\vo)$;
  \item $\Psi(u\otimes\eta_\omega) = \Psi(u)\otimes\eta_\omega$, for $u \in \tG(\FQ_\vo)$.
\end{enumerate}
\end{thm}

By \propref{prop:rep:sink}(1), the involutions $\Psi$ are well
defined. (1) follows from \propref{prop:kef:dual}. (2) follows from
\propref{prop:rep:f}(2).

\begin{thm}\label{thm:rep:basis}
The basis $\CB(\vo)$ of $\tG(\FQ_\vo)$ satisfies
\begin{enumerate}
  \item $(b,b')\in\delta_{bb'}+q^{-1}\N[q^{-1}]$ for $b,b'\in\CB(\vo)$;
  \item $\Psi(b) = b$, for $b\in\CB(\vo)$;
  \item the subset $\N[q,q^{-1}][\CB(\vo)] \subset \tG(\FQ_\vo)$ is stable under $K^{\pm1}_i, E^{(n)}_i, F^{(n)}_i$;
  \item $b\otimes\eta_\omega\in\CB(\vo\omega)$ for $b\in\CB(\vo)$.
\end{enumerate}
\end{thm}

(1) has been proved in \thmref{thm:rep:form}. (2) follows from
\propref{prop:rep:sink}(1). (3) is clear from definition. (4)
follows from \propref{prop:rep:f}(4).

\begin{rem}\label{rem:rep:main}
(1) Note that the involutions $\Psi$ from \thmref{thm:rep:bar} are
uniquely determined by the properties (1)(2) therein along with the
normalization $\Psi(\CB_\emptyset)=\CB_\emptyset$. Therefore, they
coincide with those involutions expressed in terms of
quasi-$R$-matrix in \cite{Lu93}, which share the same properties and
normalization.

(2) The basis $\CB(\vo)$ is identical to the canonical basis
introduced by Lusztig \cite[14.4.12, 27.3]{Lu93}. First, both bases
satisfy \ref{thm:rep:basis}(1)(2) and both are contained in
$\bG(\FQ_\vo)$, so they may differ at most by signs. Moreover, both
bases satisfy the recursive formula \ref{thm:rep:basis}(4).
Comparing the crystal graph on Lusztig's basis and the positivity
property \ref{thm:rep:basis}(3) of $\CB(\vo)$, one is able to see
the ambiguity of sign does not happen.

(3) The positivity result stated in \thmref{thm:rep:basis}(3) is
new. It was proved for single highest weight integrable modules in
\cite[22.1]{Lu93} when $(I,H)$ is simply-laced. It was also proved
for quantum $sl_n$ independently in \cite{Su08} by using the
categorification \cite{Su07}.

In fact, the canonical bases of the highest weight integrable
$U$-modules (they coincide with Kashiwara's global crystal bases
\cite{Ka91}) were originally specialized from a common basis, also
referred to as canonical basis or global crystal basis, of the
negative part $U^-$ of $U$ (the subalgebra generated by $F_i$, $i\in
I$). Lusztig's geometric construction of the basis of $U^-$ yields
further a positivity property \cite[Theorem 11.5]{Lu91}. From that
\thmref{thm:rep:basis}(3) is naturally expected for single highest
weight integrable modules.

However, \thmref{thm:rep:basis}(3) for tensor product modules is far
from clear from Lusztig's original algebraic construction. So our
geometric realization of canonical bases is yet another illustration
of the advantage of geometric approach to the representation theory.

(4) Nakajima \cite{Na01} associated to each module $\Lambda(\vo)$ a
Lagrangian subvariety $\tilde\FZ$ of the Nakajima's quiver variety
$\sqcup_\nu\FM_\nu$ (see also Malkin's tensor product variety
\cite{Ma03}). In fact, there is a one-to-one correspondence
\begin{equation*}
\begin{split}
  & \{ \text{``characteristic varieties'' of the simple perverse sheaves in $\FQ_\vo$} \} \\
  & \quad \quad \quad \quad \longleftrightarrow \{ \text{irreducible components of $\tilde\FZ$} \}.
\end{split}
\end{equation*}
So the category $\FQ_\vo$ can be reformulated in a more intrinsic
way as the full subcategory of $\FD$ formed by the semisimple
complexes whose characteristic varieties are contained in
$\tilde\FZ$.
\end{rem}

\subsection{Examples}

\subsubsection{Quantum $sl_2$}

The underlying finite graph consists of a single vertex:
$(I,H)=(\{i\},\emptyset)$ and $(\hI,\hH)=(\{i,\hi\},\{i\to\hi,\hi\to
i\})$. We choose $\Omega=\{i\to\hi\}$.

Assume $\nu=ri+d\hi$. Observe that $[\tE\nOi/\Gn]$ is precisely the
Grassmannian
$$Gr(r,d)=\{V\subset\BF^d\mid\dim V=r\}.$$
Moreover, the thick subcategory $\TN\nOi$ of $\D([\tE\nOi/\Gn])$ is
generated by zero. So
$$\FD_{ri+d\hi}=\D(Gr(r,d)).$$

Assume $\vo=(d_1\hi,d_2\hi,\dots,d_t\hi)$ and $d=\sum_nd_n$. Let
$P_\vo\subset GL(\BF^d)$ be the parabolic subgroup preserving the
partial flag
$$0\subset\BF^{d_1}\subset\BF^{d_1+d_2}\subset\cdots\subset\BF^{d}.$$
Note that the simple perverse sheaves in $\FQ_\vo$ are all
$P_\vo$-equivariant. Counting the number of $P_\vo$-orbits and
comparing it with the dimension of $\tG(\FQ_\vo)$ show that the
simple perverse sheaves in $\FQ_\vo$ are exactly those
$P_\vo$-equivariant simple perverse sheaves in $\D(Gr(r,d))$, $0\le
r\le d$.

Therefore, we recover the categorification described in \cite{Zh07}.

\subsubsection{Quantum $sl_3$}

The underlying finite graph consists of two vertices and a single
edge joining them. Let $\Omega,\Omega'$ be orientations as below.
\begin{equation*}
\begin{array}{cccc}
\xymatrixrowsep{1.5pc}\xymatrix{
  \\
  i \ar@<-.2ex>@_{<-}[r] \ar@<.2ex>@^{->}[r]
  & j
} \quad & \quad \xymatrixrowsep{1.5pc}\xymatrix{
  \hi \ar@<-.2ex>@_{<-}[d] \ar@<.2ex>@^{->}[d]
  & \hj \ar@<-.2ex>@_{<-}[d] \ar@<.2ex>@^{->}[d]\\
  i \ar@<-.2ex>@_{<-}[r] \ar@<.2ex>@^{->}[r]
  & j
} \quad & \quad \xymatrixrowsep{1.5pc}\xymatrix{
  \hi \ar@<-.2ex>@_{<-}[d]_{1}
  & \hj \ar@<-.2ex>@_{<-}[d]^{3} \\
  i \ar@<.2ex>@^{->}[r]^{2}
  & j
} \quad & \quad \xymatrixrowsep{1.5pc}\xymatrix{
  \hi \ar@<-.2ex>@_{<-}[d]
  & \hj \ar@<-.2ex>@_{<-}[d] \\
  i \ar@<-.2ex>@_{<-}[r]
  & j
} \vspace{1ex} \\
(I,H) \quad & \quad (\hI,\hH) \quad & \quad \Omega \quad & \quad
\Omega'
\end{array}
\end{equation*}

The whole category $\FD$ has already been too complicated to be
studied in general within an example. Below are some typical cases
for reference.

Case (1). Let $\omega=\hi+\hj$. The canonical basis of
$\Lambda(\omega)$ consists of eight elements
\begin{equation*}
\begin{array}{rrr}
  \eta_\omega, & F_i\eta_\omega, & F_j\eta_\omega, \\
  F_i^{(2)}F_j\eta_\omega, & F_j^{(2)}F_i\eta_\omega, & F_iF^{(2)}_jF_i\eta_\omega, \\
  F_jF_i\eta_\omega, & F_iF_j\eta_\omega.
\end{array}
\end{equation*}
So, totally seven $\FD_\nu$ are involved here:

(i) $\nu=\omega$. Clearly $\FD_\nu=\D([E\nO/\Gn])=\D(\pt)$. The
unique simple perverse sheaf gives rise to the basis element
$\eta_\omega$.

(ii) $\nu=i+\omega$. $[\tE\nOi/\Gn] = \BF^\times/GL(\BF) = \pt$.
Thus $\FD_\nu=\D(\pt)$. The unique simple perverse sheaf gives rise
to $F_i\eta_\omega$.

(iii) $\nu=j+\omega$. Similarly using $\Omega'$ one establishes
$\FD_\nu=\D(\pt)$. The unique simple perverse sheaf gives rise to
$F_j\eta_\omega$.

(iv) $\nu=2i+j+\omega$. $[\tE\nOi/\Gn] = [\BF/GL(\BF)]$. $\TN\nOi$
is generated by the simple perverse sheaf supported on
$[\BF^0/GL(\BF)]$. Thus
$\FD_\nu\cong\D([\BF^\times/GL(\BF)])=\D(\pt)$. The unique simple
perverse sheaf gives rise to $F_i^{(2)}F_j\eta_\omega$.

(v)  $\nu=i+2j+\omega$. Similarly $\FD_\nu\cong\D(\pt)$ and the
unique simple perverse sheaf gives rise to
$F_j^{(2)}F_i\eta_\omega$.

(vi)  $\nu=2i+2j+\omega$. The only simple perverse sheaf on
$[E\nO/\Gn]$ not contained in $\CN\nO$ is the intermediate extension
of $\bar\Q_{l,[U/\Gn]}[-1]$, where $U$ is the $\Gn$-orbit $\{ x\in
E\nO \mid \dim\Img x_1=\dim\Img x_2=\dim\Img x_3=1, \; \Ker
x_1\cap\Ker x_3x_2=0 \}$. So, $\FD_\nu$ has a unique simple perverse
sheaf, which gives rise to $F_iF_j^{(2)}F_i\eta_\omega$.

Moreover, a direct computation shows
$(F_iF^{(2)}_jF_i\eta_\omega,F_iF^{(2)}_jF_i\eta_\omega)=1$. That
being said, the unique simple perverse sheaf in $\FD_\nu$ has no
self extensions. So, $\FD_\nu\cong\D(\pt)$.

(vii)  $\nu=i+j+\omega$. We partition $E\nO$ into eight strata
$U_{d_1,d_2,d_3}=\{ x\in E\nO \mid \dim\Img x_a=d_a\}$. The simple
perverse sheaves in $\FD_\nu$ are provided by the intermediate
extensions of those on the strata $[U_{1,0,0}/\Gn]=[\pt/GL(\BF)]$,
$[U_{1,0,1}/\Gn]=\pt$, $[U_{0,1,1}/\Gn]=\pt$,
$[U_{1,1,1}/\Gn]=\BF^\times$. The intermediate extensions of
$\bar\Q_{l,[U_{1,0,1}/\Gn]}$, $\bar\Q_{l,[U_{1,1,1}/\Gn]}[1]$ give
rise respectively to $F_jF_i\eta_\omega, F_iF_j\eta_\omega$.

\smallskip

Case (2). The canonical basis of $\Lambda(\hi)\otimes\Lambda(\hj)$
is as follows.
\begin{equation*}
\begin{array}{rrr}
  (\eta_{\hi}\otimes\eta_{\hj}), & F_i(\eta_{\hi}\otimes\eta_{\hj}), & F_j(\eta_{\hi}\otimes\eta_{\hj}), \\
  F_i^{(2)}F_j(\eta_{\hi}\otimes\eta_{\hj}), & F_j^{(2)}F_i(\eta_{\hi}\otimes\eta_{\hj}), & F_iF^{(2)}_jF_i(\eta_{\hi}\otimes\eta_{\hj}), \\
  F_jF_i(\eta_{\hi}\otimes\eta_{\hj}), & F_iF_j(\eta_{\hi}\otimes\eta_{\hj}), & F_jF_i\eta_{\hi}\otimes\eta_{\hj}.
\end{array}
\end{equation*}
The first eight elements are provided by the same simple perverse
sheaves as $\Lambda(\omega)$. The last one is given by the
intermediate extension of $\bar\Q_{l,[U_{1,0,0}/\Gn]}[-1]$ in
$\FD_{i+j+\omega}$.

\smallskip

Case (3). The canonical basis of $\Lambda(\hj)\otimes\Lambda(\hi)$
is as follows.
\begin{equation*}
\begin{array}{rrr}
  (\eta_{\hj}\otimes\eta_{\hi}), & F_i(\eta_{\hj}\otimes\eta_{\hi}), & F_j(\eta_{\hj}\otimes\eta_{\hi}), \\
  F_i^{(2)}F_j(\eta_{\hj}\otimes\eta_{\hi}), & F_j^{(2)}F_i(\eta_{\hj}\otimes\eta_{\hi}), & F_iF^{(2)}_jF_i(\eta_{\hj}\otimes\eta_{\hi}), \\
  F_jF_i(\eta_{\hj}\otimes\eta_{\hi}), & F_iF_j(\eta_{\hj}\otimes\eta_{\hi}), & F_iF_j\eta_{\hj}\otimes\eta_{\hi}.
\end{array}
\end{equation*}
The first eight elements are provided by the same simple perverse
sheaves as $\Lambda(\omega)$. The last one is given by the
intermediate extension of $\bar\Q_{l,[U_{0,1,1}/\Gn]}$ in
$\FD_{i+j+\omega}$.

\subsection{Abelian categorification}\label{sec:rep:abel}

\subsubsection{}

Let $P\von$ be the set of simple perverse sheaves (up to
isomorphism) in $\FQ_\vo\cap\FD_\nu$. $P\von$ is finite since there
are only finitely many $\vm\rhd\vo$ with $|\vm|=\nu$. Let $L\von$
denote the direct sum $\bigoplus_{A\in
P\von}A\in\FQ_\vo\cap\FD_\nu$.

Define a graded $\bar\Q_l$-algebra
$$\Ao := \bigoplus_\nu\FA^\bullet\von$$
where
$$\FA^\bullet\von := \Ext^\bullet_\FD(L\von,L\von).$$

Note that for every complex $A\in\FD$,
$$\Ext^\bullet_\FD(L\von,A) \quad \text{(resp. $\Ext^\bullet_\FD(A,L\von)$)}$$
defines a graded left (resp. right) $\Ao$-module.

Let $\Ao\mof$ denote the category of finite-dimensional graded left
$\Ao$-modules and let $\Ao\pmof$ denote the full subcategory formed
by the projectives.

\subsubsection{}\label{sec:rep:abel:alg}

Basic properties of $\Ao$ and $\Ao\mof$ are as follows. The second
one follows readily from the basic properties of a t-structure. The
others are clear from the implications $(2) \Rightarrow (3)
\Rightarrow (4)(5) \Rightarrow (6)(7)$.
\begin{enumerate}
\item
$\dim \FA^\bullet\von < \infty$, which is immediate from
\propref{prop:rep:form}.
\item
$\Ao$ is $\N$-graded and $\FA^0_\vo = \bigoplus_\nu \bigoplus_{A\in
P\von} \Hom_\FD(A,A)$ each summand of which is isomorphic to
$\bar\Q_l$.
\item
The elements $\Id_A\in\FA^0_\vo$, $A\in P\von$ are precisely the
indecomposable idempotents of $\Ao$.
\item
The $\bar\Q_l$-summands of $\Ao/\FA^{>0}_\vo$ enumerate the simple
objects of $\Ao\mof$ up to grading shifts.
\item
The modules $\Ao\Id_A=\Ext^\bullet_\FD(L\von,A)$, $A\in P\von$
enumerate the indecomposable projectives of $\Ao\mof$ up to grading
shifts.
\item
The obvious map $\Hom_\FD(A,B) \to
\Hom_{\Ao}(\Ext^\bullet_\FD(L\von,A), \Ext^\bullet_\FD(L\von,B))$ is
a bijection for $A,B\in\FQ_\vo\cap\FD_\nu$.
\item
The obvious map
$\Ext^\bullet_\FD(A,L\von)\otimes_{\Ao}\Ext^\bullet_\FD(L\von,B) \to
\Ext^\bullet_\FD(A,B)$ is a bijection for
$A,B\in\FQ_\vo\cap\FD_\nu$.
\end{enumerate}

\subsubsection{}\label{sec:rep:abel:fun}

Let $\FG$ be an endofunctor of the triangulated category $\FD$.
Assume $\FG$ has a left adjoint $\FG'$ and assume both $\FG,\FG'$
restrict to endofunctors of $\FQ_\vo$. We associate to $\FG$ a
graded $\Ao$-bimodule
$$\FG^\bullet := \oplus_{\nu,\nu'} \Ext^\bullet_\FD(L\vonp,\FG L\von)
  = \oplus_{\nu,\nu'} \Ext^\bullet_\FD(\FG'L\vonp,L\von).
$$
By property \ref{sec:rep:abel:alg}(5), $\FG^\bullet$ is (left and
right) projective, therefore defines an exact endofunctor of
$\Ao\mof$ by tensoring from the left.

It follows from property \ref{sec:rep:abel:alg}(7) that
$$\FG^\bullet \otimes_{\Ao} \Ext^\bullet_\FD(L\von,A)
  = \oplus_{\nu'} \Ext^\bullet_\FD(L\vonp,\FG A), \quad
  \text{for $A\in\FQ_\vo\cap\FD_\nu$}
$$
and that for two such endofunctors $\FG_1,\FG_2$ of $\FD$,
$$\FG_1^\bullet\otimes_{\Ao}\FG_2^\bullet = (\FG_1\FG_2)^\bullet.$$

\subsubsection{}

Let $\FK^{\pm1\bullet}_i,\FE^{(n)\bullet}_i,\FF^{(n)\bullet}_i$ be
the graded projective $\Ao$-bimodules associated to the functors
$\FK^{\pm1}_i,\FE^{(n)}_i,\FF^{(n)}_i$ and regard them as exact
endofunctors of $\bG(\Ao\mof)$. This makes sense according to
\propref{prop:kef:adjoint} and \propref{prop:rep:kef}.

Then \thmref{thm:kef} implies

\begin{thm}
There are isomorphisms of graded projective $\Ao$-bimodules.
\begin{enumerate}
\setlength{\itemsep}{1.25ex}
\item
  $\FK^\bullet_i\otimes\FK^{-1\bullet}_i = \FK^{-1\bullet}_i\otimes\FK^\bullet_i = \Ao, \;\;
  \FK^\bullet_i\otimes\FK^\bullet_j = \FK^\bullet_j\otimes\FK^\bullet_i$;
\item
  $\FK^\bullet_i\otimes\FE^{(n)\bullet}_j = \FE^{(n)\bullet}_j\otimes\FK_i^{\bullet-na_{ij}}$;
\item
  $\FK^\bullet_i\otimes\FF^{(n)\bullet}_j = \FF^{(n)\bullet}_j\otimes\FK_i^{\bullet+na_{ij}}$;
\item
  $\FE^{(n-1)\bullet}_i \otimes\FE^\bullet_i \cong \bigoplus\limits_{0\le m<n} \FE_i^{(n)\bullet+n-1-2m}$;
\item
  $\FF^{(n-1)\bullet}_i \otimes\FF^\bullet_i \cong \bigoplus\limits_{0\le m<n} \FF_i^{(n)\bullet+n-1-2m}$;
\item
  $\FE^\bullet_i\otimes\FF^\bullet_i \oplus
  \bigoplus\limits_\nu \bigoplus\limits_{0\le m<\nu_i-\bar\nu_i} \FA_{\vo,\nu}^{\bullet+\nu_i-\bar\nu_i-1-2m}$
  \\
  $\cong \FF^\bullet_i\otimes\FE^\bullet_i \oplus
  \bigoplus\limits_\nu \bigoplus\limits_{0\le m<\bar\nu_i-\nu_i} \FA_{\vo,\nu}^{\bullet+\bar\nu_i-\nu_i-1-2m}$;
\item
  $\FE^\bullet_i\otimes\FF^\bullet_j \cong \FF^\bullet_j\otimes\FE^\bullet_i, \;\; i\ne j$;
\item
  $\bigoplus\limits_{\stackrel{0\le m\le 1-a_{ij}}{\text{$m$ odd}}} \FE^{(m)\bullet}_i\otimes\FE^\bullet_j\otimes\FE^{(1-a_{ij}-m)\bullet}_i
  \cong \bigoplus\limits_{\stackrel{0\le m\le 1-a_{ij}}{\text{$m$ even}}} \FE^{(m)\bullet}_i\otimes\FE^\bullet_j\otimes\FE^{(1-a_{ij}-m)\bullet}_i$;
\item
  $\bigoplus\limits_{\stackrel{0\le m\le 1-a_{ij}}{\text{$m$ odd}}} \FF^{(m)\bullet}_i\otimes\FF^\bullet_j\otimes\FF^{(1-a_{ij}-m)\bullet}_i
  \cong \bigoplus\limits_{\stackrel{0\le m\le 1-a_{ij}}{\text{$m$ even}}} \FF^{(m)\bullet}_i\otimes\FF^\bullet_j\otimes\FF^{(1-a_{ij}-m)\bullet}_i$.
\end{enumerate}
\end{thm}

\subsubsection{}

Let $\bG(\Ao\mof)$ denote the Grothendieck group of the abelian
category $\Ao\mof$. That is, $\bG(\Ao\mof)$ is the free
$\Z[q,q^{-1}]$-module defined by the generators each for an
isomorphism class of objects from $\Ao\mof$ and the relations
\begin{enumerate}
\renewcommand{\theenumi}{\roman{enumi}}
  \item $[M^\bullet]=[M'^\bullet]+[M''^\bullet]$,
  for an exact sequence $M'^\bullet \hookrightarrow M^\bullet \twoheadrightarrow M''^\bullet$;
  \item $[M^{\bullet+1}]=q^{-1}[M^\bullet]$, for $M^\bullet\in \Ao\mof$.
\end{enumerate}

\begin{cor}
The followings endow $\bG(\Ao\mof)\otimes\Q(q)$ with a structure of
$U$-module.
\begin{equation*}
\begin{split}
  & K^{\pm1}_i[M^\bullet]:=[\FK^{\pm1\bullet}_i\otimes M^\bullet], \\
  & E^{(n)}_i[M^\bullet]:=[\FE^{(n)\bullet}_i\otimes M^\bullet], \\
  & F^{(n)}_i[M^\bullet]:=[\FF^{(n)\bullet}_i\otimes M^\bullet].
\end{split}
\end{equation*}
\end{cor}

\subsubsection{}

By property \ref{sec:rep:abel:alg}(5)(6), the functor
$$\FQ_\vo \to \Ao\pmof, \quad\quad A \mapsto \oplus_\nu\Ext^\bullet_\FD(L\von,A),$$
defines an equivalence of categories, which in turn induces a linear
isomorphism
$$\bG(\FQ_\vo)\otimes\Q(q) \cong \bG(\Ao\mof)\otimes\Q(q).$$
It is clear from \ref{sec:rep:abel:fun} that this is an isomorphism
of $U$-modules.

Moreover, the above functor sends the simple perverse sheaves to the
indecomposable projectives of $\Ao\mof$
$$\bigcup_\nu \; \{ \; \Ao\Id_A=\Ext^\bullet_\FD(L\von,A) \mid A\in P\von \}.$$

Summarizing, we obtain

\begin{thm}
There is an isomorphism of $U$-modules
$$\bG(\Ao\mof)\otimes\Q(q) \cong \Lambda(\vo)$$
so that the set of indecomposable projectives of $\Ao\mof$
$$\{ \Ao e \mid \text{$e$ is an indecomposable idempotent of $\Ao$} \}$$
gives rise to the canonical basis of $\Lambda(\vo)$.
\end{thm}

\end{document}